\documentclass[11pt]{article}
 
\usepackage[utf8]{inputenc}
\usepackage{graphicx}
\usepackage{amssymb, amsmath, amsthm}
\usepackage[mathscr]{euscript}
\usepackage{pgfplots}
\usepackage{tikz}
\usetikzlibrary{shapes.misc}
\usetikzlibrary{calc}
\usetikzlibrary{matrix}
\usepackage{enumitem}

\usepackage{makeidx}  
\makeindex
\usepackage{hyperref}

\usepackage{euscript}

\usepackage{latexsym}
\usepackage{amsmath, amssymb,graphics}
\usepackage[all]{xy}
\usepackage{color}

\usepackage{tikz}
\input xy
\xyoption{all}
\usepackage{euscript}
\usepackage{multirow}
\usepackage{etex, pictexwd,dcpic}
\usetikzlibrary{positioning}
\usetikzlibrary{shapes.geometric}
\usetikzlibrary{shapes.misc}
\usetikzlibrary{calc}
\usetikzlibrary{positioning}

\usepackage{pgflibraryarrows}
\usepackage{pgflibrarysnakes}

\usetikzlibrary{trees} 
\usetikzlibrary[trees] 

\usepackage{fullpage}

\newtheorem{theorem}{Theorem}
\newtheorem{definition}[theorem]{Definition}
\newtheorem{lemma}[theorem]{Lemma}

\newtheorem{proposition}[theorem]{Proposition}

\newtheorem{remark}{Remark}

\newtheorem{example}{Example}

\numberwithin{equation}{section}

\renewcommand{\(}{\begin{equation*}}
\renewcommand{\)}{\end{equation*}}
\newcommand{\bea}{\begin{eqnarray*}}
\newcommand{\eea}{\end{eqnarray*}}
\newcommand{\R}{{\mathbb R}}

\newcommand{\Z}{{\mathbb Z}}

\def\endofproof {\hfill{$\Box$}\\}

\newcommand{\cA}{\ensuremath{\mathcal A}}

\newcommand{\cF}{\ensuremath{\mathcal F}}
\newcommand{\cG}{\ensuremath{\mathcal G}}

\newcommand{\cR}{\ensuremath{\mathcal R}}

\def\r{\rightarrow}

\newcommand{\N}{\ensuremath{\mathcal N}}
\newcommand{\beq}{\begin{equation}}
\newcommand{\eeq}{\end{equation}}

\newcommand{\onto}{\twoheadrightarrow}

\newcommand{\into}{\hookrightarrow}

\newcommand{\theproof}{\noindent {\bf Proof.\ }}

\numberwithin{equation}{section}

\renewcommand{\(}{\begin{equation}}
\renewcommand{\)}{\end{equation}}

\newcommand{\RR}{{\mathbb R}}
\newcommand{\ZZ}{{\mathbb Z}}

\def\N{{\mathbb N}}
\def\R{{\mathbb R}}
\def\Z{{\mathbb Z}}

\def\1{{\bf 1}}

\def\<{\langle}
\def\>{\rangle}

\numberwithin{equation}{section}


\newcommand{\BB}{\ensuremath{\mathbb B}}

\newcommand{\chp}{\ensuremath{\mathscr{C}\mathrm{h}^{+}}}

\newcommand{\sset}{\ensuremath{s\mathscr{S}\mathrm{et}}}
\newcommand{\set}{\ensuremath{\mathscr{S}\mathrm{et}}}

\newcommand{\sab}{\ensuremath{s\mathscr{A}\mathrm{b}}}

\newcommand{\sh}{\ensuremath{\mathscr{S}\mathrm{h}}}

\newcommand{\cartsp}{\mathscr{C}\mathrm{art}\mathscr{S}\mathrm{p}}
\newcommand{\map}{\mathrm{Map}}

\begin{document}

\title{Massey products in differential cohomology via stacks} 

 \author{
 Daniel Grady and Hisham Sati\\
  }

\maketitle

\begin{abstract} 

 We extend Massey products from cohomology to differential cohomology via stacks, 
 organizing and generalizing existing constructions in Deligne cohomology. 
 We study the properties and show how they are related to more classical 
 Massey products in de Rham, singular, and Deligne cohomology. The setting and the 
 algebraic machinery via stacks allow for computations and make the construction 
 well-suited for applications.  We illustrate with several examples from differential 
 geometry and mathematical physics.

 \end{abstract}

 \tableofcontents
 
\section{Introduction}

Massey products were introduced in \cite{Ma} and further developed and generalized 
in \cite{Kr} \cite{May}. The existence of (higher)
Massey products indicates the complexity of the topology of a
space. They also determine whether and how various characterizing properties of a 
space might be related, in particular how homotopy of a space 
might be related to its cohomology \cite{GM}. On the one hand, Massey products can be viewed
as secondary cohomology operations associated with the primary operation given by the cup product.
 On the other hand, they can also be seen 
as higher order products in homotopy ($A_\infty$) algebras (see \cite{Jim}\cite{BV}). 

\medskip
 A differential graded algebra (DGA) is a (not a
priori commutative) graded algebra $A$ with a map $d:A\r A$ of
degree $+1$ which satisfies the relations (up to sign conventions) 
$dd=0$ and 
$d(ab)= (da)b +(-1)^{\dim a}a(db)$.
Then the
cohomology $H(A)$ of $A$ with respect to $d$ is a graded algebra.
It has further certain operations called (matrix) Massey products, the 
simplest of which is a
correspondence 
\(
      H(A)\otimes H(A)\otimes H(A) \r H(A)
\) 
which is denoted by $\langle a,b,c\rangle$, where $a,b,c\in H(A)$. 
This has dimension $\dim(a)+\dim(b)+\dim(c)-1$, is 
defined only when $ab=bc=0\in H(A)$, and 
 is not well-defined but rather only defined modulo terms of
the form $ax +yb$ where $x$ and $y$ are some (auxiliary) elements of $H(A)$. 
The indeterminacy may, however, sometimes be
excluded, for example for dimension reasons, which occurs in applications. 
Generally,  we have
 $ab=dy$ and $bc=dz$ for $y,z\in A$, so that
\( \langle a,b,c \rangle= yc + (-1)^{\dim a + 1}az \) 
 is a cocycle, with the cohomology class defined
modulo the indeterminacy given above. 

\medskip
There are other notions of Massey products, but all are essentially variations 
on this principle. If the Massey product $\langle a_1,....,a_n \rangle$ exists, then all 
``lower'' Massey products necessarily vanish, although the converse is not true in general. One
may also apply a similar construction for matrices of elements, leading to matric 
Massey products \cite{May}, where notions related to formal flatness of the connection 
become important.

\medskip
Differential cohomology has played an important role recently by combining 
geometric and topological data, namely usual cohomology and differential forms, 
 in a coherent way 
   \cite{CS} \cite{EV} \cite{Bry} \cite{DF} \cite{Ga} \cite{HS} \cite{SS} \cite{Bun1} \cite{BS} \cite{Urs} \cite{BB}.
 It is natural then to 
 try to extend Massey products, which exist in both 
 of these ingredients, to differential cohomology.  
Massey products have been considered in Deligne cohomology in
 \cite{Ded} \cite{Sch} \cite{We} \cite{MT}.
We extend the definitions and constructions to the level of stacks,
\footnote{Throughout the paper, by stacks we mean simplicial sheaves 
as discussed, for instance, in \cite{Lur} \cite{DI} \cite{Urs} and recalled in 
section  \ref{Sec Diff Coh}.} 
 which has the 
virtue of allowing for 
vast generalization to various settings and to a plethora of applications. 
We believe this formulation has an advantage both for theory and for applications. 
In particular, we emphasize that desired properties and behaviour of the Massey products
are clearer and systematic in stacks, and computations are generally doable and are 
more efficient there, making them quite suitable for applications. 
The constructions are based on the thesis \cite{Dan}, but we have taken 
the opportunity here to sharpen 
the results and add properties and applications. We emphasize that this paper is 
the first part of a bigger project, aimed at developing various concrete computational 
techniques for differential cohomology theories.

\medskip
The paper is organized as follows. 
In Sec. \ref{Sec MasseyDiff}, we provide the setting for the two main
 ingredients that we would like to combine, namely classical (generalized) Massey products
 in Sec. \ref{Sec MasseyClass} and differential cohomology in Sec. \ref{Sec Diff Coh}. 
 We set up the former in the general framework of \cite{May} \cite{Kr} \cite{Bab}
 and the latter in the language of stacks (see \cite{FSSt} \cite{SSS3} \cite{Urs}).
Then we recall the Deligne-Beilinson cup product in differential cohomology,
as set forth in \cite{FSS1} \cite{FSS2}, in Sec. \ref{Sec cup}. The first encounter of 
Massey products in differential cohomology in the particular setting of Deligne cohomology
is recalled in Sec. \ref{Sec Massey hyper} and adapted slightly to our language. Our
main construction is then described in Sec. \ref{Sec Massey stacks}, where we first
set up the powerful machinery needed, in the form of the Dold-Kan 
correspondence, in Sec. \ref{doldkan}, and then provide the main definitions in
Sec. \ref{Sec Stacky Massey}. Part of this construction, together with a lot of the 
homotopic background appeared in the second author's thesis \cite{Dan}. 
A vast generalization, along the lines of the classic work of May, 
is presented in Sec. \ref{Sec General stacky}, where we present all three of differential, singular, and 
de Rham Massey products within the same setting. Then in Sec. \ref{Sec properties} we 
give the properties of the Massey products thus defined. These turn out to be 
rather attractive in general, with some unexpected features. 

\medskip
In Sec. \ref{Sec applications}, we illustrate (some aspects of)
 the construction with various applications. 
We will first give instances of where the classical Massey 
products arise in applications, and we apply our constructions in previous sections to supply the differential 
refinements of these applications.  
We extend the constructions and discussion in \cite{FSS1} from cup product Chern-Simons theories
to what we might call {\it Massey product Chern-Simons theories}. 
In Sec. \ref{Sec triv}, we illustrate how (stacky) Massey products arise in trivialization of 
higher structures, such as (differential) String, Fivebrane \cite{SSS2}, and Ninebrane 
structures \cite{9brane}. This gives natural trivializations of Chern-Simons 
theories at the level of (higher) bundles with connections.
There are two expressions that involve three differential cohomology classes, namely the 
stacky Massey
product and the triple Deligne-Beilinson cup product. A natural question is whether these
are related. Indeed, we propose such a relation via transfer in the context of cobordism. 

Then, in Sec. \ref{Sec anomaly}, we see how systems arising 
generally in anomaly cancellation lead naturally to (stacky) Massey products. 
Finally, in settings inspired by type IIA and type IIB string theory in 
Sec. \ref{Sec Bianchi} and Sec. \ref{Sec quad}, respectively, we illustrate how 
these lead to stacky Massey products. Interestingly, the latter gives rise to a 
quadruple Massey product. The reader need not be familiar with these string 
theories in order to follow the discussion.

\section{Massey products and differential cohomology}
\label{Sec MasseyDiff} 
\subsection{Classical (generalized) Massey products}
\label{Sec MasseyClass}

We recall some notions from \cite{May} \cite{Kr} \cite{Bab}. This will be useful for the applications that
we will consider later as well as a starting point for comparison with 
our stacky constructions.

\medskip
Let $({\cA},d)$ be a differential graded algebra over $\RR$
endowed with augmentation. Let $M(\cA)$ be the set of all upper
triangular half-infinite matrices with entries in $\cA$, zeroes on
the diagonal and finitely many nonzero entries, i.e.
\( 
M(\cA)=\left\{ A=(a_{ij}), a_{ij}\in \cA, a_{ij}=0 ~{\rm for}~
j\leq i~ {\rm and}~ i,j \geq n+1 ~{\rm for~ some}~ n ~\text{with}~ i, j \in \N \right\}\;. 
\label{MA}
\)
The last condition distinguishes in $M(\cA)$ a subset (which is in
fact a subalgebra) $M_n(\cA)$ consisting of all $(n\times
n)$-matrices with entries in $\cA$. The algebra $M(\cA)$ is
bigraded and endowed with a bigraded Lie bracket. We introduce the
differential $d$ on $M(\cA)$ as $dA=(da_{ij})_{i,j\geq 1}$.
The algebra $\cA$ admits an involution given by 
$a \mapsto
{\overline a}=(-1)^k a$, 
which can be extended to an
automorphism of $M(\cA)$ as ${\overline
A}=(\overline{a}_{ij})_{i,j \geq 1}$, with the differential $d$
satisfying the generalized Leibnitz rule $d(AB) = (dA)B +
{\overline A}(dB)$.
In \cite{Bab}, the Maurer-Cartan operator $\mu:~M(\cA) \longrightarrow M(\cA)$ was defined as
$\mu(A)=dA - {\overline A} \cdot A$.
Then a matrix $A\in M(\cA)$ is said to be a matrix of formal connection
if it satisfies the Maurer-Cartan equation in $\cA$,
\( 
dA -{\overline A} \cdot A \equiv 0 ~~{\rm mod} \ker A\;, 
\) i.e.
$A$ is a formal connection if $\mu(A) \in \ker A$. Here $\ker A$
is a $\cA$-module generated by matrices ${1_{ij}}$ such that $A\cdot
1_{ij}=1_{ij} \cdot A$, where $1_{ij}$ denotes the matrix that has all zero entries
except for 1 as the $ij$-entry.  
Note that this implies that $AB=BA$ for any matrix $B \in
\ker A$. The element  $\mu(A)$ is called the curvature of the formal
connection $A$, and can be shown to be closed (see e.g. \cite{Bab} \cite{Ch}).

\medskip
Now comes the relation between Maurer-Cartan and Massey
products. The generalized Massey products are the cohomology
classes of the curvature matrices of the formal connection $A$,
i.e. if $A$ is a solution to the Maurer-Cartan equation then the
entries of the matrix $[\mu(A)]$ are the generalized Massey
products \cite{Bab}. Geometrically, this means that the latter
measure the deviation of connections from flat ones, so that the
connection is flat if they vanish. Later we will make use of this approach 
in describing Massey products in stacks. 

\medskip
Classical Massey products in integral cohomology $H^*(X; \Z)$ arise 
by taking  ${\cal A}$ to be an algebra over the commutative ring $\Z$, 
with the multiplication being associative but not necessarily 
graded-commutative. 
Now let $\alpha, \beta, \gamma$ be the cohomology classes of
closed elements $a\in \cA^p$, $b \in \cA^q$, and $c\in \cA^r$. The triple Massey product
$\langle \alpha, \beta, \gamma \rangle$ is defined if one can solve
the Maurer-Cartan equation with the formal connection
\begin{displaymath}
A= \left( \begin{array}{cccc}
0 & a & {\tilde f} & \ast \\
0 & 0 & b & {\tilde g} \\
0 & 0 & 0 & c \\
0 & 0 & 0 & 0
\end{array} \right)\;.
\end{displaymath}
This is equivalent to the two separate equations 
\(
d{\tilde f} = (-1)^p a \wedge b
\quad {\rm and} \quad 
d{\tilde g} = (-1)^q a \wedge c \label{fg}
\)
and that implies that the Massey product is defined if and only if
\( 
\alpha \cup \beta = \beta \cup \gamma=0 \in H^*(\cA)\;. 
\)
The matrix $\mu(A)$ has the form
\begin{displaymath}
\mu(A)=dA- {\overline A} \cdot A= \left( \begin{array}{cccc}
0 & 0 & 0 & \tau \\
0 & 0 & 0 & 0 \\
0 & 0 & 0 & 0 \\
0 & 0 & 0 & 0
\end{array} \right)
\end{displaymath}
and defines the Massey product $[\mu(A)]$ which is equal to the cohomology class
\(
\langle \alpha, \beta, \gamma \rangle= [\tau]= \left[ (-1)^{p+1} a
\wedge {\tilde g} + (-1)^{p+q}{\tilde f} \wedge c \right]\;. 
\)
Here $[a] \in H^*({\cal A})$ denotes 
the cohomology class of a closed element $a \in {\cal A}$, and 
$[A]=([a_{ij}])_{i, j \geq1}\in M(H^*({\cal A})$, for a closed matrix $A \in M({\cal A})$,
 denotes the corresponding matrix whose entries are the cohomology classes
 of the entries $a_{ij}$ of $A$. 
Since ${\tilde f}$ and ${\tilde g}$ are defined by expressions \eqref{fg} up
to closed elements from $\cA$, the triple Massey product $\langle
\alpha, \beta, \gamma \rangle$ is defined modulo $\alpha\cdot
H^{q+r}(\cA) + \gamma \cdot H^{p+q}(\cA)$.

\subsection{Differential cohomology}
\label{Sec Diff Coh} 


There are several different approaches to differential cohomology. Initially, 
we will be concerned with the construction as Deligne cohomology \cite{Bry} \cite{Ga}. 
We will then move to the stacky setting, which illuminates the true nature of differential 
cohomology as a theory which counts  isomorphism classes of higher $U(1)$-gerbes with connection
(generalizing the usual discussion for the gerbe case in \cite{Bry}). 

\medskip
The classical construction relies on hypercohomology of a complex of objects of an abelian category as
 an extension to complexes of the usual cohomology of an object.
For $n\in \N$, let $\Z^\infty_{\mathcal{D}}[n]$ be the sheaf of chain complexes given by
$$
\Z^\infty_{\cal D}[n]:=[\hdots \to 0 \to \underline{\ZZ} \into \Omega^{0} \to \Omega^{1} \to \hdots \to \Omega^{n-1}]\;,
$$
where $\underline{\ZZ}$ is in degree
\footnote{This is a descending grading, which is the opposite of the usual grading of the de Rham complex. 
That is, we are viewing this as a chain complex rather than a cochain complex. Furthermore, we take $V[n]$ 
to denote the chain complex shifted by $n$, so that $V$ is in degree $n$.}
 $n$ and $\Omega^{n-1}$ is the sheaf of real-valued $(n-1)$-forms in degree 
$0$. Given a manifold $X$, the degree $n$ sheaf hypercohomology with coefficients in $\Z^\infty_{\cal D}[n]$ can be 
 defined to be the degree $n$ differential cohomology of $X$:
\(
\widehat{H}^{n}(X;\ZZ):=H^{n}(X; \Z^\infty_{\mathcal{D}}[n])\;.
\)
If $X$ is paracompact, then these cohomology groups are given by the cohomology of the total complex 
of the {\v Cech}-Deligne double complex corresponding to a good open cover of $X$. In what follows, we 
will always assume that $X$ is paracompact, so that the hypercohomology groups can be computed by 
either taking arbitrary injective resolutions, or via this more explicit {\v C}ech approach. 

\medskip
In \cite{SS} (see also \cite{Bun1}), it was observed that these cohomology groups fit nicely into 
an exact hexagon
\(\label{diamond}
\begin{tikzpicture}[baseline=(current bounding box.center)]
\matrix(m) [matrix of math nodes, column sep=3em, row sep=3em, nodes={anchor=center}]{
 & \Omega^{n-1}(X)/\mathrm{im}(d) & & \Omega^{n}_{\mathrm{cl}}(X) &
\\
H^{n-1}_{\mathrm{dR}}(X) & & \widehat{H}^{n}(X;\ZZ) & & H^{n}_{\rm dR}(X)\;,
\\
 & H^{n-1}(X;U(1)) & & H^{n}(X,\ZZ) &
\\
};

\path[->]
(m-1-2) edge node[above] {\small $d$} (m-1-4)
(m-2-1) edge (m-1-2)
(m-1-2) edge node[above] {\small $a$} (m-2-3)
(m-2-3) edge node[above] {\footnotesize $I$} (m-3-4)
(m-2-1) edge (m-3-2)
(m-3-2) edge (m-3-4)
(m-3-4) edge (m-2-5)
(m-3-2) edge (m-2-3)
(m-2-3) edge node[right] {\ \footnotesize $R$} (m-1-4)
(m-1-4) edge (m-2-5);
\end{tikzpicture}
\)
where the bottom row is the Bockstein sequence and the diagonals are exact. The map $R$ is called the \emph{curvature} 
map and $I$ is called the \emph{integration} map. Notice that, by exactness, in the case that the curvature of a differential 
cohomology class vanishes, the class lies in the image of the inclusion $H^{n-1}(X;U(1))\into \widehat{H}^{n}(X;\ZZ)$. 
We call these classes \emph{flat}, as they represent $n$-gerbes with connections of vanishing curvature. Differential 
cohomology therefore detects the topological information -- when the class is flat -- and the differential geometric 
information encoded by the curvature. See \cite{CS} \cite{EV} \cite{Bry} \cite{DF} \cite{Ga} \cite{HS} \cite{SS} 
\cite{BS} \cite{Urs} \cite{BB}  for more details on the various approaches. 

\medskip
As we mentioned earlier, our point of view henceforth will be mainly that of stacks. 
We will recall and introduce some stacks that will be useful for us. 
We start by surveying some basic concepts and definitions from \cite{Lur} \cite{FSSt} \cite{FSS1} \cite{FSS2}, 
adapted to our setting. We will provide only as much detail as necessary to introduce our stacks. 

\medskip
For $n \in \N$, let $\cartsp$ be the category with objects convex open subsets of 
Cartesian space $\RR^{n}$ (hence diffeomorphic to $\R^n$), and morphisms smooth functions. 
A \emph{smooth prestack} is simply a functor
$$
\cF:\cartsp^{\rm op}\to \sset
$$
with target the category of simplicial sets. The passage from prestacks to stacks is achieved by imposing a 
sort of gluing condition on $\cF$. Roughly speaking, a stack $\cF$ attaches an entire space (equivalently 
simplicial set) of data to each object in $\cartsp$. This data should be viewed as being \emph{local} data. 
The gluing condition then assembles this data into a geometric object, which is a stack. 
More precisely, we say that a prestack $\cF$ satisfies \emph{descent} if for each $U\in \cartsp$ 
and each open cover $\{U_{i}\}_{i\in I}$ of $U$ with contractible finite intersections 
$U_{i_1i_2\hdots i_k}$, we have a weak equivalence
\(
\label{desc}
\cF(U)\simeq {\rm holim}\Big\{\begin{tikzpicture}
[descr/.style={fill=white,inner sep=2.5pt}, back line/.style={densely dotted},
cross line/.style={preaction={draw=white, -,
line width=6pt}}, baseline={(current bounding box.center)}]
\matrix (m) [matrix of math nodes,row sep=3em, column sep=2em,
text height=1.5ex,text depth=0.25ex]{
 \hdots  & \prod_{i,j,k}\cF(U_{ijk}) & \prod_{i,j}\cF(U_{ij}) & \prod_{i}\cF(U_{i})  \\
};

\path[-stealth]
([yshift=1.5pt]m-1-2.west) edge ([yshift=1.5pt]m-1-1.east)
([yshift=-1.5pt]m-1-2.west) edge ([yshift=-1.5pt]m-1-1.east)
([yshift=4.5pt]m-1-2.west) edge ([yshift=4.5pt]m-1-1.east)
([yshift=-4.5pt]m-1-2.west) edge ([yshift=-4.5pt]m-1-1.east)
([yshift=3pt]m-1-3.west) edge ([yshift=3pt]m-1-2.east)
([yshift=-3pt]m-1-3.west) edge ([yshift=-3pt]m-1-2.east)
([yshift=0pt]m-1-3.west) edge ([yshift=0pt]m-1-2.east)
([yshift=1.5pt]m-1-4.west) edge ([yshift=1.5pt]m-1-3.east)
([yshift=-1.5pt]m-1-4.west) edge ([yshift=-1.5pt]m-1-3.east);
\end{tikzpicture}\Big\}\;.
\) 
In particular, if $\cF$ takes values in Kan complexes, this weak equivalence is part of 
an actual \emph{homotopy} equivalence. The reader may notice the following.

\begin{itemize}
\item If we change the target category to $\set$ and impose the stronger condition that the \emph{strict} limit over the diagram was \emph{isomorphic} to $\cF(U)$, we would recover the gluing condition for a sheaf. 
\item If we change the target category to groupoids, then the above condition recovers the usual notion of descent for classical stacks. 
\end{itemize}
In the latter, homotopy equivalence is simply categorical equivalence of groupoids. Hence the gluing condition respects the correct notion of equivalence (which is weaker than isomorphism). We can therefore view the equivalence \eqref{desc} as the more general gluing condition for $\infty$-groupoids (or Kan complexes). We need the following (see \cite{Lur} \cite{DI} \cite{Urs}). 

\begin{definition}
We call a smooth prestack $\cF$ a \emph{smooth stack} if it satisfies descent. We denote the full subcategory 
of smooth stacks by
$$
\sh_{\infty}(\cartsp)\into [\cartsp^{\rm op},\sset]\;,
$$
where the brackets denote the category of contravariant functors from $\cartsp$ to $\sset$, with morphisms that are natural 
transformations. 
\end{definition}

Note that the above functor category is simplicially enriched in a natural way. 
Observe that for objects $X$ and $Y$ in any (locally small) category, ${\rm hom}(X, Y)$ is always a set.
This allows us to form the mapping {\it space} (i.e. simplicial set), which at level $n$ is 
$$
\big(\text{Map}(X, Y) \big)_n:=\text{hom}\big(X \times \underline{\Delta[n]}, Y\big)\;,
$$
when $X$ is fibrant and $Y$ cofibrant (this requires a model structure). Here the operation $\times$ is the Cartesian 
product in stacks, and the underline on $\Delta[n]$ denotes taking the locally constant stack 
associated to $\Delta[n]$.

\begin{remark} The inclusion functor admits a left adjoint $L$ which preserves homotopy colimits (in fact, a left Quillen adjoint \cite{Urs}). We call this functor $L$ the \emph{stackification} functor and call the image of a prestack $\cF$ under $L$ the \emph{stackification} of $\cF$.
\end{remark}

In \cite{FSSt}, the moduli stack of $n$-gerbes with connection, $\BB^{n} U(1)_{\rm conn}$, was introduced. This stack was obtained as the stackification of the $n$-prestack obtained by applying the Dold-Kan map (see section \ref{doldkan}) to the Deligne presheaf of chain complexes
$$
\ZZ^{\infty}_{\mathcal{D}}[n+1]:=[... \to 0 \to \underline{\ZZ} \into \Omega^{1} \to \Omega^{2} \to ... \to \Omega^{n}]\;.
$$
These stacks are the differential analogues of Eilenberg-MacLane spaces and, for a fixed manifold $X$, there is a bijective correspondence (a ``representation") 
\(\label{homotopycat}
\widehat{H}_{\mathcal{D}}^{n+1}(X;\ZZ)\simeq \pi_0\map(X,\BB^{n} U(1)_{\rm conn})\;,
\)
where the right hand side is the set of morphisms in the homotopy category of stacks.

\begin{remark} 
In general the right hand side of the correspondence \eqref{homotopycat} may not be well-defined. 
In order to be able to take homotopy groups of the mapping space $Y:=\map(X,\BB^{n} U(1)_{\rm conn})$,
$Y$ has to be a Kan complex, which is the case when 
$X$ is confibrant and $\BB^{n} U(1)_{\rm conn}$ 
is fibrant. 
However, since 
$\BB^{n}U(1)_{\rm conn}$ satisfies descent, it is fibrant in a particular local model structure on presheaves (see \cite{FSSt}).
Even though $X$ can be viewed as a stack, it is not cofibrant, and so we need to cofibrantly replace it.
 Indeed, if $X$ is a (paracompact) manifold, thought of as a smooth stack, with good open cover 
 $\{U_{i}\}_{i\in I}$, then we can replace $X$ by its {\v C}ech nerve
\(
C(\{U_i\}):={\rm hocolim}\Big\{\begin{tikzpicture}
[descr/.style={fill=white,inner sep=2.5pt}, back line/.style={densely dotted},
cross line/.style={preaction={draw=white, -,
line width=6pt}}, baseline={(current bounding box.center)}]
\matrix (m) [matrix of math nodes,row sep=3em, column sep=2em,
text height=1.5ex,text depth=0.25ex]{
 \coprod_{i}U_{i} & \coprod_{i,j}U_{ij} & \coprod_{i,j,k}U_{ijk} & \hdots \\
};

\path[-stealth]
([yshift=1.5pt]m-1-4.west) edge ([yshift=1.5pt]m-1-3.east)
([yshift=-1.5pt]m-1-4.west) edge ([yshift=-1.5pt]m-1-3.east)
([yshift=4.5pt]m-1-4.west) edge ([yshift=4.5pt]m-1-3.east)
([yshift=-4.5pt]m-1-4.west) edge ([yshift=-4.5pt]m-1-3.east)
([yshift=3pt]m-1-3.west) edge ([yshift=3pt]m-1-2.east)
([yshift=-3pt]m-1-3.west) edge ([yshift=-3pt]m-1-2.east)
([yshift=0pt]m-1-3.west) edge ([yshift=0pt]m-1-2.east)
([yshift=1.5pt]m-1-2.west) edge ([yshift=1.5pt]m-1-1.east)
([yshift=-1.5pt]m-1-2.west) edge ([yshift=-1.5pt]m-1-1.east);
\end{tikzpicture}\Big\}
\) 
which is both cofibrant and weak equivalent to $X$ in the category of smooth stacks ${\rm Sh}_\infty (\cartsp)$ \cite{Dug}.
 For purely model category theoretic reasons it then follows that 
$\map(C(\{U_{i}\}),\BB^{n}U(1)_{\rm conn})$
is a Kan complex and we can take $\pi_{0}$, obtaining the set of morphisms in the homotopy category. This motivates the definition
$$\map(X,\BB^{n}U(1)_{\rm conn}):=\map(C(\{U_{i}\}),\BB^{n}U(1)_{\rm conn})\;.$$
\end{remark}

\medskip
As explained in \cite{FSSt}, these stacks also have a nice geometric interpretation. The following example illustrates the point quite well.

\begin{example} 
Let $X$ be a manifold. Let us calculate the set of vertices of the mapping space $\map(X,\BB^{2}U(1)_{\rm conn})$. Using the pointwise formula for the homotopy colimit \cite{FSSt}, we have
\begin{eqnarray}
\hom(X,\BB^{2}U(1)_{\rm conn}) &=& \hom(C(\{U_{i}\}),\BB^2U(1)_{\rm conn})
\nonumber \\
 &=& \hom\Big(\int^{k\in \Delta}\underline{\Delta[k]}\times \coprod_{\alpha_{1},..,\alpha_{k}}U_{\alpha_{1},..,\alpha_{k}},~\BB^{2}U(1)_{\rm conn}\Big)
\nonumber \\
&=& \int_{k\in \Delta}\prod_{\alpha_{1},..,\alpha_{k}}\hom(\underline{\Delta[k]}\times U_{\alpha_{1},..,\alpha_{k}},~\BB^{2}U(1)_{\rm conn})
\nonumber\\
&=& \prod_{\alpha_{1},..,\alpha_{k}}\int_{k\in \Delta}\hom\big(\underline{\Delta[k]},~\BB^{2}U(1)_{\rm conn}(U_{\alpha_{1},..,\alpha_{k}})\big)\;.
\end{eqnarray}
An element of the hom in the last line can be written out explicitly as a choice maps
{\small
\begin{align}
B_{\alpha} &:\Delta[0]\to \prod_{\alpha}\BB^{2}U(1)(U_{\alpha})
\nonumber \\
A_{\alpha\beta} &: \Delta[1]\to \prod_{\alpha\beta}\BB^{2}U(1)(U_{\alpha\beta})
\nonumber \\
g_{\alpha\beta\gamma} &: \Delta[2] \to \prod_{\alpha\beta\gamma}\BB^{2}U(1)(U_{\alpha\beta\gamma})\;,
\end{align}
}
such that the face inclusions of each map are equal to their corresponding restrictions to higher intersections. Now since equivalent stacks will produce the same cohomology groups, we do not distinguish between equivalent stacks. In particular, using the exponential quasi-isomorphism, we could have equivalently defined $\BB^{2}U(1)_{\rm conn}$ to be the stackification of the prestack given by applying the Dold-Kan functor to the presheaf of chain complexes
$$
[ 0 \to \hdots \to C^{\infty}(-,U(1)) \overset{d\log}{\longrightarrow} \Omega^{1} \to \Omega^{2}]\;.
$$ 
We can therefore describe the choices of $B_{\alpha}$,$A_{\alpha\beta}$ and $g_{{\alpha\beta\gamma}}$ via 
the 2-simplex
$$
{\small 
\begin{tikzpicture}[scale=1.25]

\coordinate [label=left:$B_{{\alpha}}$] (A) at (-1.0cm,-1.0cm);
\coordinate [label=right:$B_{\gamma}$] (B) at (1.0cm,-1.0cm);
\coordinate [label=above:$B_{\delta}$] (C) at (0cm,1.0cm);
\coordinate [label=below:$g_{\alpha\beta\gamma}$] (D) at (0cm, 0cm);
\draw (A) -- node[below] {$A_{\alpha\beta}$} (B) -- node[right] {$A_{\gamma\delta}$} (C) -- node[left] {$A_{\delta\alpha}$} (A);

\end{tikzpicture}
}
$$
Here, $g_{\alpha\beta\gamma}$ is a choice of smooth $U(1)$-valued function on triple
 intersections, $A_{\alpha\beta}$ is a choice of 1-form on double intersections and 
 $B_{\alpha}$ is a choice of 2-form on open sets. Moreover, we have that these assignments
 must satisfy the conditions 
\begin{enumerate}[label= {\bf (\roman*)}]
\item $g_{\alpha\beta}~g_{\gamma\beta}^{-1}~g_{\gamma\alpha}=1$;

\item $g_{\alpha\beta\gamma}^{-1}~dg_{\alpha\beta\gamma}=d\log(g)_{\alpha\beta\gamma}=A_{\alpha\beta}-A_{\gamma\beta}+A_{\gamma\alpha}$;

\item $B_{\beta}-B_{\alpha}=dA_{\alpha\beta}$.
\end{enumerate}
We identify this data as precisely giving a gerbe with connection \cite{Bry}. Moreover, the fact that $\BB^{n}U(1)_{\rm conn}$ 
is a stack ensures that $F_{\alpha}=dB_{\alpha}$ is a globally defined 3-form: the curvature of the gerbe.
Notice that these are only the vertices in the mapping space. The entire mapping space keeps track of more information, 
namely the homotopies and higher homotopies between gerbes. These encode automorphisms in the sense of gauge
 transformations (see \cite{FSS1} \cite{FSS2}).
\end{example}

\begin{example} 
Let $X$ be a paracompact manifold and $C(\{U_i\})$ the {\v C}ech nerve of some good open cover. The maps
$$
L:C(\{U_i\})\to \BB U(1)_{\rm conn}
$$
are in bijective correspondence with circle bundles on $X$ equipped with a connection. In fact, using the calculations in the above 
example shows that such a morphism gives the data $U(1)$-valued functions $g_{\alpha\beta}$ on intersections satisfying
$
g_{\alpha\beta}g_{\beta\gamma}^{-1}g_{\gamma\delta}=1
$
on triple intersections, along with 1-forms $A_{\alpha}$ on open sets satisfying
$
A_{\alpha}-A_{\beta}=d\log(g)_{\alpha\beta}
$
on double intersections. 
If the homotopy class of $L$ is trivial, then the circle bundle is \emph{trivializable}. In fact, the trivializing map $\phi$ is nothing but a homotopy $\phi:L \to 0$. To identify this homotopy, we use the Dold-Kan correspondence. In particular, an edge in 
$\map(C(\{U_i\}),\BB U(1)_{\rm conn})$ is, by adjunction, an edge in the simplicial set
\(
\map(C(\{U_i\}),\BB^{n}U(1)_{\rm conn})={\rm DK}(\hom_{\chp}(N(C(\{U_i\})),\ZZ^{\infty}_{\mathcal{D}}[2]))\;,
\)
where $N$ is the normalized Moore functor. 
Recall that this functor gives an equivalence of categories, from simplicial abelian groups 
$\sab$ to chain complexes in non-negative degrees ${\rm Ch}_\bullet^+$
(see \cite{GJ}).
The hom in positively graded chain complexes is the truncated total complex of the {\v C}ech-Deligne double complex
$$
[\hdots \to \mathrm{tot}^{1}C(\mathcal{U},\ZZ^{\infty}_{\mathcal{D}}[2])
\to
 Z\left(\mathrm{tot}^{2}C(\mathcal{U},\ZZ^{\infty}_{\mathcal{D}}[2])\right)]\;,
$$
where $Z$ denotes the group of cocycles in that degree. Recalling that the differential is given by $D:=d+(-1)^{k}\delta$, where 
$\delta$ takes the alternating sum of restrictions, we identify an edge connecting $L$ and $0$ as an assignment of {\v C}ech-Deligne cochain $h$ of degree 1 such that $(d-\delta)h=L$. Explicitly, 
this means a choice of $U(1)$-valued function $h_{\alpha}$ on open sets such that
\begin{enumerate}[label= {\bf (\roman*)}]
\item $h_{\alpha}h_{\beta}^{-1}=g_{\alpha\beta}$;
\item $-ih_{\alpha}^{-1}dh_{\alpha}=d\log(h_{\alpha})=A_{\alpha}$.
\end{enumerate}
 A straightforward calculation shows that the pattern continues and that null homotopies of $n$-gerbes 
 (equivalently $n$-bundles, equivalently maps into $\BB^{n}U(1)_{\rm conn}$) can again be identified with trivializations.
\end{example}

Motivated by this last example, we will often refer to null homotopies as \emph{trivializations}. To summarize, the mapping space 
$\map(X,\BB^{n}U(1)_{\rm conn})$ can be identified with the set of all $n$-gerbes with connection, along with isomorphisms between these and higher homotopies between these isomorphisms. 

\begin{remark}\label{stacks}
There are several other stacks  related to $\BB^{n}U(1)_{\rm conn}$ which are useful for us and are defined as follows (see \cite{FSSt},\cite{FSS1},\cite{FSS2},\cite{Urs}):

\noindent {\bf (i)}  If we forget about the connection on the these $n$-bundles, we obtain the bare moduli stack of $n$-gerbes $\BB^{n}U(1)$. Explicitly, this stack is obtained by applying the Dold-Kan functor to the sheaf of chain complexes 
$C^{\infty}(-,U(1))[n]$: the sheaf of smooth $U(1)$-valued functions in degree $n$.

\item {\bf (ii)} We also define a stack which represents \emph{flat} $n$-bundles with connection, $\flat\BB^{n}U(1)$. This stack is obtained by applying Dold-Kan to the sheaf of chain complexes ${\rm disc}{U(1)}[n]$: the sheaf of \emph{locally constant} $U(1)$ valued functions in degree $n$. 
\footnote{
Here ``disc" refers to the underlying discrete topology. As an operation on stacks, ${\rm disc}$ is the composite functor 
$\xymatrix{{\rm disc}: {\rm Sh}_\infty \ar[r]^{{\rm ev}_*} & \sset \ar[r]^{\underline{( \cdot )}} & {\rm Sh}_\infty}$, where
${\rm ev}_*$ is the evaluation at a point and $\underline{(\cdot )}$ takes the locally constant stack associated to a simplicial set. For a smooth manifold $X$, the resulting stack ${\rm disc}(X)$ is sometimes denoted instead by $\underline{X}^{\delta}$.
}

\item {\bf (iii)} We have a stack representing the truncated de Rham complex $\flat_{\rm dR}\BB^{n}U(1)$ obtained by applying Dold-Kan to the truncated de Rham sheaf of chain complexes 
$$
\Omega^{\leq n}_{\rm cl}:=[\hdots 0 \to \Omega^{0}\to \Omega^{1} \to \hdots \to \Omega^{n}_{\rm cl}]\;.
$$

\item {\bf (iv)} Finally, we define the stack of closed $n$-forms $\Omega^{n}_{\rm cl}$ to be the stack obtained by applying Dold-Kan to the sheaf of closed $n$-forms.
\end{remark}

One way to see that the second stack really does detect \emph{flat} $n$-gerbes with connection is to observe that, by Poincar\'e lemma, one has a quasi-isomorphism of sheaves
$$
{\rm disc}(U(1))[n]\simeq [ 0 \to \hdots \to C^{\infty}(-,U(1)) \overset{d\log}{\longrightarrow} \Omega^{1} \to \hdots \to \Omega^{n}_{\rm cl}]\;,
$$
where on the right we have \emph{closed} $n$-forms in degree $0$. These $n$-forms are to be interpreted as giving the connection on the corresponding bundle. Hence, if the form is closed then the bundle is flat.

\medskip
The moduli stack $\BB^{n}U(1)_{\rm conn}$ is related to the stacks in Remark \ref{stacks} in various ways.
 In \cite{FSSt}\cite{Urs}, it was observed that $\BB^{n}U(1)_{\rm conn}$ is the homotopy pullback
\(\label{curv}
\xymatrix{
\BB^{n}U(1)_{\rm conn} \ar[d] \ar[rr]^-R && \Omega^{n+1}_{\rm cl} \ar[d]^\iota 
\\
\BB^{n}U(1) \ar[rr]^-\theta && \flat_{\rm dR}\BB^{n+1}U(1)\;,
}
\)
where the left composite $\BB^{n}U(1)_{\rm conn}\to \BB^{n}U(1)\overset{\theta}{\to} \flat_{\rm dR}\BB^{n+1}U(1)$ is homotopic to the map
\(
{\rm curv}:\BB^{n}U(1)_{\rm conn}\to \flat_{\rm dR}\BB^{n+1}U(1)
\)
induced by the morphism of sheaves of chain complexes
\(
 \xymatrix{
 \ZZ \ar[rr]^-i \ar[d]^-i && \Omega^{0} \ar[rr] \ar[d]^d && \cdots \ar[d]^d \ar[rr] && \Omega^n \ar[d]^d \\
\Omega^{0} \ar[rr]^-d && \Omega^{1} \ar[rr] && \cdots \ar[rr] && \Omega^{n+1}_{\rm cl}\;.
}
 \)
This map gives the full de Rham data for the curvature of a bundle with connection. In fact, if one calculates
 the sheaf hypercohomology in degree $0$ of the bottom row, say via the {\v C}ech-de Rham 
 complex (as in \cite{BT}), one gets  $H^{n}_{\rm dR}(X)$. Consequently, the map ${\rm curv}$ induces a map
\(
{\rm curv}_*:\pi_{0}\map(X,\BB^{n}U(1)_{\rm conn})\longrightarrow 
H^{n+1}_{\rm dR}(X)\;,
\)
which sends an $(n-1)$-gerbe with connection to the de Rham class of its curvature. The following 
proposition might certainly be known to experts, 
but we include a proof for completeness.
\begin{lemma}\label{flatfiber} The homotopy fiber of the map
$$
R:\BB^{n}U(1)_{\rm conn}\onto \Omega^{n+1}_{\rm cl}
$$
can be identified with $\flat \BB^{n}U(1)$.
\end{lemma}
\theproof
The map $R$ is induced by the morphism of sheaves of chain complexes
\(
 \xymatrix{
 \ZZ \ar[rr]^-i \ar[d]^-0 && \Omega^{0} \ar[rr] \ar[d]^0 && \cdots \ar[d]^0 \ar[rr] && \Omega^{n-1} \ar[d]^d \\
{0} \ar[rr]^-d && 0 \ar[rr] && \cdots \ar[rr] && \Omega^{n}_{\rm cl}\;.
}
 \)
 Since this map is degree-wise surjective by Poincar\'e lemma (traditionally in highest form-degree, and trivially in lower degrees), 
it is a fibration in the projective model structure on presheaves of chain complexes. We can therefore calculate the homotopy fiber as the kernel of that map. By inspection, the kernel is
$$
[\hdots \to \underline{\Z} \into \Omega^{0}\to \Omega^{1} \to \hdots \to \Omega^{n}_{\rm cl}]\;,
$$
which, via the exponential map, is quasi-isomorphic to 
$$
[\hdots C^{\infty}(-,U(1))\overset{d {\rm log}}{\longrightarrow} \Omega^{1}\to \hdots \to \Omega^{n}_{\rm cl}]\;.
$$
Again, by Poincar\'e lemma, this sheaf of chain complex is quasi-isomorphic to ${\rm disc}(U(1))[n].$
Since the Dold-Kan functor is a right Quillen adjoint and preserves weak equivalences, it takes fibration sequences 
to fibration sequences and we have the desired result.
\endofproof

Using the above proposition along with diagram \eqref{curv} and the pasting lemma for homotopy pullbacks, 
we observe that we have the following iteration of homotopy pullbacks \cite{Urs}
\(
\label{diffhex}
\xymatrix{
{\flat \BB^{n-1}U(1) } \ar[r] \ar[d] & {\BB^{n}U(1)} \ar[r] \ar[d] & \ast \ar[d]^-0
&
\\
\ast \ar[r]^-0 \ar[d] & \flat_{\rm dR}\BB^{n-1}U(1) \ar[r] \ar[d] & \flat \BB^{n}U(1) \ar[r] \ar[d] & \ast \ar[d]^-0
\\
\ast \ar[r]^-0 & \Omega^{\leq n-1} \ar[r] \ar[d] & \BB^{n}U(1)_{\rm conn} \ar@{->>}[r] \ar[d] & \Omega^{n+1}_{\rm cl}  \ar[d]
\\
 & \ast \ar[r]^-0 &  \BB^{n}U(1) \ar[r] & \flat_{\rm dR} \BB^{n+1}U(1)\;,
}
\)
where $0$ is the $0$ map. 
From Lemma \ref{flatfiber} along with this last diagram, we immediately get the following:
\begin{proposition}
The based loop stack $\Omega\BB^{n}U(1)_{\rm conn}$ can be identified with the stack $\flat \BB^{n-1}U(1)$.
\label{Prop loop}
\end{proposition}
\theproof
Consider the homotopy pullback square
$$
\xymatrix{
 \flat \BB^{n-1}U(1) \ar[r] \ar[d] & \BB^{n}U(1) \ar[r] \ar[d] & \ast \ar[d]^-0
\\
\ast \ar[r]^-0 \ar[d] & \flat_{\rm dR}\BB^{n-1}U(1) \ar[r] \ar[d] & \flat \BB^{n}U(1) \ar[d] 
\\
\ast \ar[r]^-0 & \Omega^{\leq n-1} \ar[r] & \BB^{n}U(1)_{\rm conn} 
}
$$
within diagram \eqref{diffhex}. Such a homotopy pullback, given by the outer square, can be taken as a definition of the loop space. Alternatively, a homotopy pullback can be computed explicitly as the paths in $\BB^nU(1)_{\rm conn}$ connecting 
the point inclusion $\ast\to \BB^nU(1)_{\rm conn}$ to itself: a loop.
\endofproof

Note that Massey products in the homology of the based loop space is classically 
considered in \cite{St} \cite{Ch2}. The above discussions allows us to recast the 
``differential cohomology diamond" using our stacks.  

\begin{proposition} The differential cohomology diagram \eqref{diamond} lifts to a diagram of stacks
\(
\label{stackdiamond}
\begin{tikzpicture}
\matrix (m) [matrix of math nodes, column sep=3em, row sep=3em, nodes={anchor=center}]{
& \Omega^{\leq n-1} & & \Omega^{n}_{\rm cl}&
\\
 \flat_{\rm dR}\BB^{n-1}U(1) & &  \BB^{n}U(1)_{\rm conn} & & \flat_{\rm dR} \BB^{n}U(1) 
\\
& \flat \BB^{n}U(1) & & \BB^{n}U(1) &
\\
};

\path[->]
(m-1-2) edge node[above] {\small $d$} (m-1-4)
(m-2-1) edge (m-1-2)
(m-1-2) edge node[above] {\small$a$} (m-2-3)
(m-2-3) edge[->>] node[above] {\footnotesize $I$} (m-3-4)
(m-2-1) edge (m-3-2)
(m-3-2) edge node[above] {\footnotesize $\beta$} (m-3-4)
(m-3-4) edge (m-2-5)
(m-3-2) edge node[above] {\small $j$} (m-2-3)
(m-2-3) edge[->>] node[right] {~ \footnotesize $R$} (m-1-4)
(m-1-4) edge (m-2-5);     
;
\end{tikzpicture}
\)
where the diagonals are fibration sequences.
\end{proposition}
\theproof
This is the same diagram as a portion of diagram \eqref{diffhex} rotated.
The top and bottom horizontal maps in \eqref{stackdiamond} are defined as
the compositions $d=Ra$ and $\beta=jI$. Fixing a manifold $X$, mapping into 
this diagram, and passing to connected components, i.e. taking 
$\pi_0 \map (X, -)$, we recover the diamond diagram \eqref{diamond}. 
Note that $d$ in \eqref{stackdiamond} recovers the usual exterior
derivative, by the nature of $R$, and that $\beta$ recovers the Beckstein 
by uniqueness of the latter as a cohomology operation. 
\endofproof

We now explain how to go the other direction, i.e. from stacks to Deligne cohomology. 
We have seen that for a manifold $X$, the mapping space $\map(X,\BB^{n}U(1)_{\rm conn})$ can be identified with the space 
of $n$-gerbes equipped with connections (along with all isomorphisms and higher isomorphisms between them). It will be convenient to
 organize this mapping space itself into a stack. We define the \emph{mapping stack} to be the stackification of the prestack 
 given by the assignment
\(
U\mapsto \map(X\times U,\BB^nU(1)_{\rm conn})
\)
for each $U\in \cartsp$. We denote this stack by 
\footnote{Note that this is not to be confused with homotopy classes of maps as the notation might suggest.}
$[X,\BB^nU(1)_{\rm conn}]$.

\begin{remark}
Notice the following:

\noindent {\bf (i)} If we evaluate the mapping stack on the terminal object in $\cartsp$ (the point) and take $\pi_0$, we recover the usual differential cohomology groups from the correspondence \eqref{homotopycat}
$$
\pi_0[X,\BB^nU(1)_{\rm conn}](*)\simeq \pi_0\map(X\times *,\BB^nU(1)_{\rm conn})\simeq \widehat{H}^n(X,\ZZ)\;.
$$

\item {\bf (ii)} Since the mapping stack is clearly functorial in both arguments and the stackification functor preserves homotopy fibers (it is left exact), we can map into the diagram \eqref{diffhex} to obtain the diagram
$$
\hspace{-2mm}
\begin{tikzpicture}
\matrix (m) [matrix of math nodes, column sep=2em, row sep=3em, nodes={anchor=center}]{
& \left[X, \Omega^{\leq n-1}\right] & & \left[X, \Omega^{n}_{cl}\right]&
\\
 \left[X, \flat_{\rm dR}\BB^{n-1}U(1)\right] & &  \left[X, \BB^{n}U(1)_{\rm conn}\right] & & \left[X, \flat_{\rm dR} \BB^{n}U(1)\right] 
\\
& \left[X, \flat \BB^{n}U(1)\right] & & \left[X, \BB^{n}U(1)\right] &
\\
};

\path[->]
(m-1-2) edge node[above] {\small $d$} (m-1-4)
(m-2-1) edge (m-1-2)
(m-1-2) edge node[above] {\small $a$} (m-2-3)
(m-2-3) edge[->>] node[above] {\footnotesize $I$} (m-3-4)
(m-2-1) edge (m-3-2)
(m-3-2) edge (m-3-4)
(m-3-4) edge (m-2-5)
(m-3-2) edge (m-2-3)
(m-2-3) edge[->>] node[right] {~ \footnotesize $R$} (m-1-4)
(m-1-4) edge (m-2-5);     
;
\end{tikzpicture}
$$
where the diagonals are again fibration sequences. If we evaluate this previous diagram at the point and apply $\pi_0$, we indeed reproduce the usual differential cohomology diamond diagram \eqref{diamond}.
\end{remark}

\subsection{Cup product in differential cohomology}
\label{Sec cup} 


Deligne \cite{De} and Beilinson \cite{Be} showed that differential cohomology admits a distinguished cup product 
refining the usual cup product on singular cohomology. This product is defined on sections of $\Z^\infty_{\cal D}[n]$ 
by the formula
\(\label{dbcup}
\alpha \cup_{\rm DB} \beta=\left \{ 
\begin{array}{ccc}
\alpha \beta, && \ \mathrm{deg}(\alpha)=n
\\
\alpha \wedge d\beta, && \ \mathrm{deg}(\alpha)=0
\\
0, && \ \text{otherwise}.
\end{array}\right.
\)   
Note that the grading here is such that the first case is simply multiplication by an integer. In fact, it is obvious from 
the definition that the Deligne-Beilinson (henceforth DB) cup product composed with the natural inclusion
$$\underline{\ZZ}[n]\into \ZZ^{\infty}_{\mathcal{D}}[n]$$
simply multiplies the two locally constant integer-valued functions. Since the sheaf cohomology of the locally constant 
sheaf $\underline{\ZZ}$, equipped with this product, is simply the ordinary cohomology ring with integral coefficients, 
one immediately sees that this cup product does indeed refine the usual cup product. 

\medskip
Equipped with this cup product, $\widehat{H}^{*}(X;\ZZ)$ becomes an associative and graded-commutative
 ring \cite{Bry}. This cup product structure also refines the wedge product of forms in the sense that the curvature map 
 $R:\widehat{H}^{*}(X;\ZZ)\to \Omega^{*}_{\rm cl}$ defines a homomorphism of graded commutative rings \cite{Bun1}. 
 In particular this implies that the cup product of two classes of odd degree is flat. It can also be shown \cite{Bun1} that 
 the cup product of a flat class with any other class is again flat and that the inclusion of $H^*(X, U(1))$ into 
 $\widehat{H}^*(X; \Z)$ is a two sided ideal. 

\medskip
We now turn to the cup product, viewed as a morphism of stacks. In \cite{FSS1} it was observed that the lax monoidal 
structure of the Dold-Kan map gives rise to a cup product, exhibited as a morphism
\(
\cup:\BB^{m} U(1)_{\rm conn} \times \BB^{n} U(1)_{\rm conn}\longrightarrow \BB^{n+m+1} U(1)_{\rm conn}
\)
of stacks. This map is obtained by simply taking the DB cup product \eqref{dbcup}
$$
\cup_{\rm DB}:\ZZ^{\infty}_{\mathcal{D}}[n+1]\otimes \ZZ^{\infty}_{\mathcal{D}}[m+1] \longrightarrow
 \ZZ^{\infty}_{\mathcal{D}}[n+m+2]\;,
$$
applying the Dold-Kan map
$$
{\rm DK}(\cup_{\rm DB}): {\rm DK}(\ZZ^{\infty}_{\mathcal{D}}[n+1]\otimes \ZZ^{\infty}_{\mathcal{D}}[m+1]) \longrightarrow
 {\rm DK}(\ZZ^{\infty}_{\mathcal{D}}[n+m+2])\;,
$$
and using the lax monoidal structure $\varphi$ of the map ${\rm DK}$ to get a map
$$
\cup= {\rm DK}(\cup_{\rm DB})\circ \varphi: {\rm DK}(\ZZ^{\infty}_{\mathcal{D}}[n+1])\times {\rm DK}(\ZZ^{\infty}_{\mathcal{D}}[n+1]) \to
 {\rm DK}(\ZZ^{\infty}_{\mathcal{D}}[n+1]\otimes \ZZ^{\infty}_{\mathcal{D}}[m+1]) \to {\rm DK}(\ZZ^{\infty}_{\mathcal{D}}[n+m+2])\;.
$$
Applying the stackification functor then gives the desired map. This map then induces a map of stacks
(which we also denote as $\cup$)
\(
\cup:[X,\BB^{n}U(1)_{\rm conn}]\times [X,\BB^{m}U(1)_{\rm conn}]
\longrightarrow
 [X,\BB^{n+m+1}U(1)_{\rm conn}]\;.
\)
The following two propositions are implicit in \cite{FSS1} \cite{FSS2}. 

\begin{proposition} 
\label{BD refines}
The DB cup product refines the singular cup product. That is, we have a commutative diagram
$$
\begin{tikzpicture}
\matrix (m) [matrix of math nodes, column sep=3em, row sep=3em]{
\BB^{n}U(1)_{\rm conn}\times \BB^{m}U(1)_{\rm conn} & \BB^{n+m+1}U(1)_{\rm conn}
\\
\BB^{n+1}\ZZ \times \BB^{m+1}\ZZ & \BB^{n+m+2}\ZZ\;.
\\
};

\path[->]
(m-1-1) edge node[above] {\footnotesize $\cup_{\rm DB}$} (m-1-2)
(m-1-1) edge node[right] {\footnotesize $I\times I$} (m-2-1)
(m-2-1) edge node[above] {\footnotesize $\cup$} (m-2-2)
(m-1-2) edge node[right] {\footnotesize $I$} (m-2-2);
\end{tikzpicture}
$$
\end{proposition}
\theproof Let $p:\ZZ^\infty_{\cal D}[n+1]\to \underline{\ZZ}[n+1]$ be the projection map
$$
\xymatrix{
\underline{\Z} \ar[d]^-{\rm id} \ar[rr]^-i  && \Omega^{0} \ar[rr] \ar[d]^-0 
&& \cdots \ar[rr] \ar[d]^0 && \Omega^n \ar[d]^-0 
\\
\underline{\Z} \ar[rr]^-d  && 0 \ar[rr] && \cdots \ar[rr]  && 0 \;.
}
$$
Then, by definition of the DB cup product, the diagram
$$
\begin{tikzpicture}
\matrix (m) [matrix of math nodes, column sep=3em, row sep=3em]{
\ZZ^{\infty}_{\mathcal{D}}[n+1]\otimes \ZZ^{\infty}_{\mathcal{D}}[m+1] & \ZZ^{\infty}_{\mathcal{D}}[n+m+2]
\\
\underline{\ZZ}[n+1]\otimes \underline{\ZZ}[m+1] & \underline{\ZZ}[n+m+2]
\\
};

\path[->]
(m-1-1) edge node[above] {\footnotesize $\cup_{\rm DB}$} (m-1-2)
(m-1-1) edge node[right] {\footnotesize $p$} (m-2-1)
(m-2-1) edge node[above] {\footnotesize $\cup$} (m-2-2)
(m-1-2) edge node[right] {\footnotesize $p$} (m-2-2);
\end{tikzpicture}
$$
commutes in sheaves of chain complexes. Applying the Dold-Kan functor and using naturality of the lax monoidal 
structure map gives the result.
\endofproof
\begin{proposition} 
\label{DB refines2}
The cup product refines the wedge product, and we have a commutative diagram
$$
\begin{tikzpicture}
\matrix (m) [matrix of math nodes, column sep=3em, row sep=3em]{
\BB^{n}U(1)_{\rm conn}\times \BB^{m}U(1)_{\rm conn} & \BB^{n+m+1}U(1)_{\rm conn}
\\
\Omega^{n+1}_{\rm cl} \times \Omega^{m+1}_{\rm cl} & \Omega^{n+m+2}_{\rm cl}\;.
\\
};

\path[->]
(m-1-1) edge node[above] {\footnotesize $\cup_{\rm DB}$} (m-1-2)
(m-1-1) edge node[right] {\footnotesize $R\times R$} (m-2-1)
(m-2-1) edge node[above] {\footnotesize $\wedge$} (m-2-2)
(m-1-2) edge node[right] {\footnotesize $R$} (m-2-2);
\end{tikzpicture}
$$
\end{proposition}
\theproof
Let $\alpha$ and $\beta$ be sections of $\ZZ^{\infty}_{\mathcal{D}}[n+1]$ and $\ZZ^{\infty}_{\mathcal{D}}[m+1]$, respectively. 
Applying the curvature $R$ to the DB cup product \eqref{dbcup} gives
$$
R(\alpha\cup_{\rm DB}\beta)=\left\{
\begin{array}{cc}
\alpha d(\beta) & ~~\text{if}~~ \mathrm{deg}(\alpha)=n
\\
d(\alpha) \wedge d(\beta) & ~~~\text{if}~~ \mathrm{deg}(\beta)=0
\\
0 & ~~~\text{otherwise},
\end{array}\right.
$$
which is $R(\alpha)\wedge R(\beta)$. We therefore have a commuting diagram
$$
\begin{tikzpicture}
\matrix (m) [matrix of math nodes, column sep=3em, row sep=3em]{
\ZZ^{\infty}_{\mathcal{D}}[n+1]\otimes \ZZ^{\infty}_{\mathcal{D}}[m+1] 
& \ZZ^{\infty}_{\mathcal{D}}[n+m+2]
\\
\Omega_{\rm cl}^{n+1}\otimes \Omega_{\rm cl}^{m+1} & \Omega^{n+m+2}_{\rm cl}\;.
\\
};

\path[->]
(m-1-1) edge node[above] {\footnotesize $\cup_{\rm DB}$} (m-1-2)
(m-1-1) edge node[right] {\footnotesize $R$} (m-2-1)
(m-2-1) edge node[above] {\footnotesize $\wedge$} (m-2-2)
(m-1-2) edge node[right] {\footnotesize $R$} (m-2-2);
\end{tikzpicture}
$$
Applying the Dold-Kan map DK gives the result in stacks.
\endofproof

The above results show that, 
in general, the Deligne-Beilinson cup product does not refine the de Rham wedge product 
for the whole de Rham complex, but does so only for the top and bottom degrees. 
However, for the triple product  the only cup products that arise
are between degree zero and degree one cocycles, so that nothing is missed in passing to $\cup_{DB}$. 
We will make this more precise in Prop.  \ref{Prop Msing}.

\subsection{Massey products in hypercohomology}

\label{Sec Massey hyper}


Massey products in Deligne-Beilinson cohomology are described in 
\cite{Ded} \cite{Sch} \cite{We} \cite{MT}. In this section, we review the construction
 for hypercohomology found in \cite{Sch}, with a slightly adapted language for later 
 comparison and generalization. 
In section \ref{Sec Massey stacks} 
we generalize this construction in two ways, which we describe. We use the 
Dold-Kan correspondence to establish these products 
in the stacky setting. We also use the machinery of May \cite{May} to exhibit these products as differential 
matric Massey products.  

\medskip
Let $R$ be a commutative ring and let $\mathcal{C}^{\bullet}(n)$, $n\in \N$, be a sequence of positively graded chain complexes 
of $R$-modules. Moreover, let us assume that this sequence comes equipped with maps
$$
\cup:\mathcal{C}^{\bullet}(n)\otimes \mathcal{C}^{\bullet}(m)\to \mathcal{C}^{\bullet}(n+m)\;,
$$
which are associative in the sense that 
\(
\label{ass}
\cup \circ (id \otimes \cup)=\cup \circ ( \cup \otimes id)\;.
\)
The maps $\cup$ induce an associative product on cohomology 
$$
\cup: H^{\bullet}(n)\otimes H^{\bullet}(m)\to H^{\bullet}(n+m)\;,
$$
called the \emph{cup product}. Once a well-defined notion of a
cup product is established, one can define the Massey products via
the following.

\begin{definition} 
Let $l\geq 2$ and let $n_{1}, \cdots ,n_{l}$ and
  $m_{1}, \cdots, m_{l}$ be integers. Define
$$
n_{s,t}=\sum_{i=s}^{t}(n_{i}-1)\ ~~\text{and}~~ \
m_{s,t}=\sum_{i=s}^{t}m_{i} , \quad {\rm for} ~~1\leq s \leq t \leq l\;,
$$ 
and let $\bar{a}=(-1)^{q+1}a$ denote the
twist of a class $a\in \mathcal{C}^{q}(n)$. We define the $l$-fold Massey
product as follows:

\vspace{0.1in}

\noindent {\bf (i).} Let $a_{i}\in H^{m_{i}}(\mathcal{C}^{\bullet}(n_{i}))$ be
cohomology classes. Suppose there exists cochains $a_{s,t}\in
\mathcal{C}^{m_{s,t+1}}(n_{s,t})$ such that $a_{i,i}$ is a representative of
$a_{i}$ and that
$$d a_{s,t}=\sum_{i=s}^{t-1}\bar{a}_{s,i}\cup a_{i+1,t} \quad {\rm for} \quad 1\leq s\leq t \leq l,~~ (s,t)\neq (1,l)\;.$$ 
We call the collection
$\mathcal{M}=\{a_{s,t}\}$ a \emph{defining system} for the $l$-fold Massey
product.

\vspace{0.1in}

\noindent {\bf (ii).} The cochain
$$a_{1,l}:=\sum_{i=1}^{l-1}\bar{a}_{1,i}\cup a_{i+1,l}\in \mathcal{C}^{m_{1,l+2}}(n_{1,l})$$
is a cocycle and represents a cohomology class $m_{l}$. We call this class the $l$-fold
\emph{Massey product} of the elements $a_{1},..,a_{l}$ with defining
system $\mathcal{M}$.

\end{definition}

In general, we would like to eliminate the dependance of the product on the
defining system. The case of $l=3$ will be the most important for
us, and in this case we are indeed able to eliminate this dependence.
The following three examples are known, and we record them to highlight
how Massey products arise in the different settings that we consider,
and how stacks will provide, in a sense, a unifying theme. 
Note that, while the above construction is fairly general,
it is not obvious how to generalize to other settings and how to do computations 
easily with it, and that is why we later use the stacky perspective.

\begin{example}

Let $a_{1},a_{2}$ and $a_{3}$ be cohomology classes as above. Suppose we have a defining system $\mathcal{M}=\{a_{s,t}\}$. This means, by
definition, that we have the relations
$$
d a_{1,2}=\bar{a}_{1,1}\cup a_{2,2} \qquad 
{\rm and}  \qquad d a_{2,3}=\bar{a}_{2,2}\cup a_{3,3}\;.
$$
Now a class $m_{l}$ representing the Massey product of this defining system has as a representing cocycle
$$a_{1,l}=\bar{a}_{1,1}\cup a_{2,3}+\bar{a}_{1,2}\cup a_{3,3}\;.$$ 
Notice that, in this case, the
class $m_{l}$ only depends on the defining system up to
cocycles. That is, for another defining system
$\mathcal{N}=\{b_{s,t}\}$, the classes $a_{2,3}-b_{2,3}$ and
$a_{1,2}-b_{1,2}$ are cocycles. Moreover, if these cocycles are coboundaries, then the Massey products of both defining systems agree. We can therefore define a Massey product, not depending on the defining system, as the quotient

$$\langle a_{1},a_{2},a_{3}\rangle \in \frac{H^{m_{1}+m_{2}+m_{3}-1}(\mathcal{C}^{\bullet}(n_{1}+n_{2}+n_{3}))}{H^{m_{1}+m_{2}-1}(\mathcal{C}^{\bullet}(n_{1}+n_{2}))\cup
  a_{3}+a_{1} \cup H^{m_{2}+m_{3}-1}(\mathcal{C}^{\bullet}(n_{2}+n_{3}))}\;.$$ 
\end{example}


\medskip
\begin{example}
Let $X$ be a smooth manifold and let $C(n)=\Omega^{*}(X)$ for each $n$, where $\Omega^{*}(X)$ is the algebra of differential forms on $X$. Let $a,b$ and $c$ be de Rham cohomology classes of degree $p$, $q$, $r$ respectively, such that $a\wedge b=0=b\wedge c$. Choose representing closed forms
$\alpha$, $\beta$, $\gamma$ for $a$, $b$, $c$ respectively, and let
$\eta$ and $\rho$ be cochains such that
$$
d\eta= \alpha \wedge \beta \quad {\rm and} \quad
d\rho=\beta \wedge \gamma\;.
$$
Then the combination 
$$
\eta \wedge \gamma -(-)^p \alpha\wedge \rho
$$
is a closed form representing the triple Massey product of $a$, $b$ and $c$ corresponding to the defining system 
$\mathcal{M}=(\alpha, \beta, \gamma, \rho, \eta)$.
 Eliminating the dependence on $\mathcal{M}$ gives a well-defined class in the quotient group
$$
H_{\rm dR}^{p+q+r-1}(X)/\left(a \cup H_{\rm dR}^{q+r-1}(X) + c\cup H_{\rm dR}^{p+q-1}(X)
\right)\;.
$$
\end{example}

\medskip
The following constitutes our initial transition to differential cohomology, which we will
develop in stacks in the following section.

\begin{example}
Consider the Deligne complex given by the sheaf of chain complexes
$$
\ZZ^{\infty}_{\mathcal{D}}[n]:=[\underline{\ZZ} \into \Omega^{0} \to \Omega^{1} \to ... \to \Omega^{n-1}]\;.
$$
 Let $X$ be a paracompact manifold with good open cover $\{U_i\}_{i\in I}$ and let 
 $C(n):=\mathrm{tot}C^{\bullet}(\{U_i\},\mathcal{\ZZ}^{\infty}_{\mathcal{D}}[n])$ be the total complex of 
 the {\v C}ech-Deligne double complex. The degree $n$ cohomology of this total complex calculates the 
 differential cohomology of $X$:
$$
\hat{H}^{n}(X;\ZZ)=H^{n}\big(\mathrm{tot}C^{\bullet}(\mathcal{U},\mathcal{\ZZ}^{\infty}_{\mathcal{D}}[n])\big)\;.
$$
The Deligne-Beilinson cup product is defined as a morphism
$$
\cup_{\rm DB}:\ZZ^{\infty}_{\mathcal{D}}[n]\otimes \ZZ^{\infty}_{\mathcal{D}}[m]\longrightarrow
 \ZZ^{\infty}_{\mathcal{D}}[n+m+1]\;,
$$
which on sections is given by the formula \eqref{dbcup}.
This map induces cup product morphisms on the total complexes $C(n)$ which are associative in the sense of the 
identity in \eqref{ass}. 
\end{example}

We can therefore use this cup product to define the Massey product in differential cohomology, viewed as the sheaf hypercohomology of the Deligne complex. Since our point of view will subsume this construction, we will delay explicit examples until sections 
\ref{Sec Massey stacks} and \ref{Sec applications}.

\section{Massey products in the language of higher stacks}
 \label{Sec Massey stacks}


We provide our main construction of stacky Massey products in this section. We start with 
setting up the machinery needed.

\subsection{The Dold-Kan correspondence} 
\label{doldkan}


The Dold-Kan correspondence will be an important component in defining the Massey product in stacks. We will use the correspondence to organize the homotopies involved in certain homotopy commuting diagrams in an algebraic way.


\medskip
The classical Dold-Kan correspondence describes an equivalence of categories (see e.g. \cite{GJ})
\begin{equation}
\begin{tikzpicture}[baseline=(current bounding box.center)]
\matrix (m) [matrix of math nodes, column sep=3em, row sep=3em, nodes={anchor=center}]{
\Gamma:\chp & \sab:N
\\
};

\path[->]
([yshift=0.1cm]m-1-1.east) edge ([yshift=0.1cm]m-1-2.west)
([yshift=-0.1cm]m-1-2.west) edge ([yshift=-0.1cm]m-1-1.east);
\end{tikzpicture}
\end{equation}
between positively graded chain complexes and simplicial abelian groups. By post-composing 
with the free-forgetful adjunction, one obtains an adjunction 
\begin{equation}
\begin{tikzpicture}[baseline=(current bounding box.center)]
\matrix (m) [matrix of math nodes, column sep=3em, row sep=3em, nodes={anchor=center}]{
{\rm DK}:=U\Gamma:\chp & \sset:NF\;.
\\
};

\path[->]
([yshift=0.1cm]m-1-1.east) edge ([yshift=0.1cm]m-1-2.west)
([yshift=-0.1cm]m-1-2.west) edge ([yshift=-0.1cm]m-1-1.east);
\end{tikzpicture}
\end{equation}
In fact, one can say more. This adjunction is a Quillen adjunction of model categories, with the projective model structure on chain complexes and the Quillen model structure on simplicial sets. As such, it preserves the homotopy theories in both categories; it  therefore comes as no surprise that for a positively graded chain complex $C_{\bullet}$ one has an isomorphism
\(
H_{n}C_{\bullet}\simeq \pi_{n}{\rm DK}(C_{\bullet})\;.
\)
For convenience, we remind the reader what the functor ${\rm DK}$ does to a chain complex, as this will be a frequently used tool in producing abelian stacks.

\medskip
Let $\Delta$ denote the category of linearly ordered sets of $n$ elements with order preserving maps. Let $C_{\bullet}$ be a positively graded chain complex. The degree $n$ component of the simplicial abelian group ${\rm DK}(C_{\bullet})$ is given by
$$
{\rm DK}(C_{\bullet})_{n}=\bigoplus_{[n]\onto[k]}C_{k}\;.
$$
Here the indexing set is taken to be all \emph{surjections} $[n]\onto [k]$. It is a bit trickier to describe the face and degeneracy maps. Let $d^{i}:[n-1]\into [n]$ be a coface map in $\Delta$. We want to define the corresponding face map. To get a map out of the direct sum, it suffices to describe the map on each factor. Therefore, we need only define the face map on a term $C_{k}$ given by a surjection $\sigma:[n]\onto [k]$. To see where to send this term, we form the composite $\sigma d^{i}[n-1]\into [n] \onto [k]$. Now this morphism need not be surjective, so we factorize  $\mu\sigma^{\prime}[n-1]\onto [m] \into [k]$ where the first map is a surjection and the second map is an injection. Then $\sigma^{\prime}$ corresponds to a term $C_{m}\into \bigoplus_{[n-1]\to[m]}A_{m}={\rm DK}(C_{\bullet})_{n-1}$. We send the factor $C_{k}$ to the factor $C_{m}$ by a map $\mu^{\prime}:C_{k}\to C_{m}$. This map is given by
\(
\mu^{\prime}=\left\{ \begin{array}{cc}
{\rm id}, &\ \mu={\rm id},
\\
(-1)^{k}d, &\ \mu=d^{k},
\\
0, &\ \text{otherwise.}
\end{array}\right.
\)
A similar construction is used to define the degeneracy maps. The following example illustrates the point quite well.

\begin{example}
Consider the chain complex $A[1]$,  with the abelian group $A$ in degree $1$ and $0$'s in all other degrees. We want to compute 
${\rm DK}(A[1])$. Using the above formula, we see that the only nonzero terms in degree $n$ are given by the surjections $[n]\onto [1]$. Each surjection can be thought of as being given by an element $i\in [n]$ which divides the set into two subsets: those that go to $0$ and those that go to $1$. We therefore have $n$ surjections and
$$
{\rm DK}(A[1])_{n}=\bigoplus_{i=1}^{n}A\;.
$$
For a coface map $d^{j}:[n-1]\to [n]$, the corresponding face map $d^{j}$ is given as follows. Let $A_{i}$ denote the copy of $A$ corresponding to the $i$th surjection. Then 
$$d_{j}(A_{i})=\left\{\begin{array}{cc}
    A_{i-1} & \text{if}\ \ i>j\neq 0, n
\\
    A_{i} & \text{if}\ \ i\leq j\neq 0, n
\end{array}\right.\;, \quad
d_{0}(A_{i})=\left\{\begin{array}{cc}
    A_{i-1} & \text{if}\ \ i\neq 0
\\
    0 & \text{if}\ \ i=0
\end{array}\right.\;, \quad
d_{n}(A_{i})=\left\{\begin{array}{cc}
    A_{i} & \text{if}\ \ i\neq n
\\
    0 & \text{if}\ \ i=n.
\end{array}\right.$$
Notice that for $j\neq 0, n$, the term corresponding to  $i=j$ and $i=j+1$ both go to the same copy of $A$. We therefore have a map $A\times A \to A$ extending the identity on each component. Hence, this morphism is just group multiplication. From this, one can see that this simplicial abelian group is just the delooping group $BA$.
\end{example}
Another way to describe the simplicial set ${\rm DK}(C_{\bullet})$, which is perhaps more conceptual, is via a labeling of simplices with elements of the chain complex $C^{\bullet}$. A 2-simplex in ${\rm DK}(C^{\bullet})$, for example, is a simplex with face, edges and vertices labeled by elements of $C^{\bullet}$
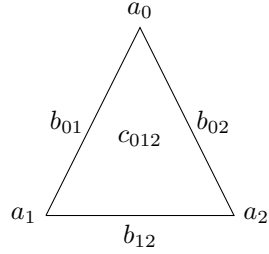
$$
{\small 
\begin{tikzpicture}[scale=1.25]

\coordinate [label=left:$a_{1}$] (A) at (-1.0cm,-1.0cm);
\coordinate [label=right:$a_{2}$] (B) at (1.0cm,-1.0cm);
\coordinate [label=above:$a_{0}$] (C) at (0cm,1.0cm);
\coordinate [label=below:$c_{012}$] (D) at (0cm, 0cm);
\draw (A) -- node[left] {$b_{01}$} (C) -- node[right] {$b_{02}$} (B) -- node[below] {$b_{12}$} (A);

\end{tikzpicture}
}
$$

\vspace{-2mm}
\noindent such that
$$
dc_{012}=b_{01}+b_{12}-b_{02} \qquad \text{and} \qquad db_{ij}=a_{j}-a_{i}\;.
$$
Here $d$ is the chain complex differential. Notice that a 2-simplex in ${\rm DK}(C_{\bullet})$, defined as before, can be
 identified as such a labeled simplex. To see this, let us calculate the data involved in specifying a 2-simplex. 
First, observe that there is exactly one surjection $0:[2]\onto [0]$, ${\rm id}:[2]\onto [2]$, and exactly two surjections $\sigma_{i}:[2]\onto [1]$. Therefore, a 2-simplex is given by a quadruple $(a_{0},b_{01},b_{02},c_{012})$, where $a_{0}$ is in degree 1, while 
$b_{01},b_{02}$ are in degree 2 (corresponding to $\sigma_{1},\sigma_{2}$, respectively), and $c_{012}$ is in degree 3. 
To determine the edges, we evaluate $d_{i}$ on this quadruple. For $i=0$, we have the following epi-mono factorizations
$$
{\rm id}\circ d^{0} = d^{0}\circ {\rm id}, \qquad \quad \sigma_{1}\circ d^{0} = d^{0}\circ 0, 
\qquad \quad \sigma_{2}\circ d^{0} = {\rm id}\circ {\rm id}\;.
$$
It follows from the formula, that the 0 face is $b_{02}$. For $i=1$, we have 
$$
{\rm id}\circ d^{1} = d^{1}\circ {\rm id},
\qquad \quad \sigma_{1}\circ d^{1} = {\rm id}\circ {\rm id},
 \qquad \quad \sigma_{2}\circ d^{1} = {\rm id}\circ {\rm id}
$$
and the 1 face is $b_{01}+b_{02}$. Finally, for $i=2$, we have
$$
{\rm id}\circ d^{2} = d^{2}\circ {\rm id},
 \qquad \quad \sigma_{1}\circ d^{2} = {\rm id}\circ {\rm id}, 
\qquad \quad \sigma_{2}\circ d^{2} = d^{1}\circ 0
$$
and the 2 edge is $dc_{012}+b_{01}$. Forming the boundary of the simplex, we get
$$
\partial (a_{0},b_{01},b_{02},c_{012})=b_{02}-(b_{01}+b_{02})+(dc_{012}+b_{01})=dc_{012}\;.
$$
That the edges of the simplex satisfy the second condition above is a straightforward
calculation and will be omitted. In fact, it is a straightforward calculation to show that 
the boundary of a general $n$-simplex must 
be equal to $d$ applied to the labeling on its $n$-face.
 
\begin{remark} 
This second description provides a powerful conceptual advantage; namely, that the differential of the 
chain complex can be viewed as obstructing the chain from being a cycle. For example, if the resulting 
simplicial set were the nerve of a groupoid, then all simplices for $n\geq 2$ would be cycles.
\end{remark}

\subsection{Stacky Massey products}
\label{Sec Stacky Massey}


We are now ready to define Massey products in the category of stacks. We begin with a discussion on 
Massey triple products and then generalize to $l$-fold Massey products. 

\medskip
The Massey triple product can be viewed as a homotopy built out of the associativity diagram 
of the cup product of three elements. In fact, suppose one is given a triple of higher gerbes with 
connection on a manifold $X$. These gerbes are given by the data 
$\mathcal{G}_{i}:X\to \BB^{n_{i}}U(1)$, $i=1,2,3$. Suppose, moreover, that these gerbes are 
chosen so that the cup products $\mathcal{G}_{1}\cup \mathcal{G}_{2}$ and $\mathcal{G}_{2}\cup \mathcal{G}_{3}$ 
are homotopic to $0$, with trivializing homotopies $\phi_{1,2}$ and $\phi_{2,3}$. In this case, we can 
build a \emph{loop} trivializing the triple product. To see this, consider the associativity diagram for the cup product

$$
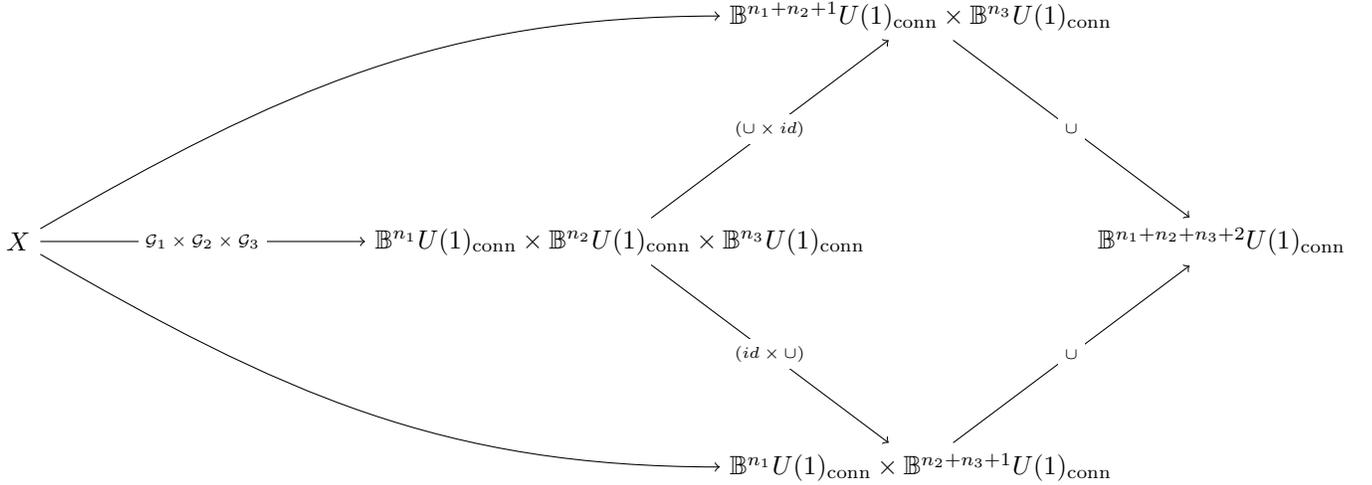
\hspace{-5mm}
\begin{tikzpicture}[descr/.style={fill=white,inner sep=2.5pt}]
\node (1) at (0,0) {\small $\BB^{n_{1}}U(1)_{\rm conn}\times \BB^{n_{2}}U(1)_{\rm conn}\times \BB^{n_{3}}U(1)_{\rm conn}$}; 
\node (2) at (4,3) {\small $\BB^{n_{1}+n_{2}+1}U(1)_{\rm conn}\times \BB^{n_{3}}U(1)_{\rm conn}$};
\node (3) at (4,-3) {\small $\BB^{n_{1}}U(1)_{\rm conn}\times \BB^{n_{2}+n_{3}+1}U(1)_{\rm conn}$};
\node (4) at (8,0) {\small $\BB^{n_{1}+n_{2}+n_{3}+2}U(1)_{\rm conn}$};
\node (5) at (-8,0) {\small $X$};

\draw[->] (1) to node[descr] {\tiny $(\cup\times id)$} (2);
\draw[->] (1) to node[descr] {\tiny $(id \times \cup)$}(3);
\draw[->] (2) to node[descr] {\tiny $\cup $}(4);
\draw[->] (3) to node[descr] {\tiny $\cup $}(4);
\draw[->] (5) to node[descr] {\tiny $\mathcal{G}_{1}\times \mathcal{G}_{2} \times \mathcal{G}_{3}$} (1);
\draw[->] (5) to[out=-30,in=180] (3);
\draw[->] (5) to[out=30, in=180] (2);
\end{tikzpicture} 
$$

\vspace{5mm}
Although the outer two maps agree, there is still nontrivial homotopy theoretic information contained in the diagram. To see this, suppose $\mathcal{G}_{1}\cup\mathcal{G}_{2}$ and $\mathcal{G}_{2}\cup\mathcal{G}_{3}$ are trivializable with trivializations $\phi_{1,2}$ and $\phi_{2,3}$, respectively. Then we can add these homotopies to the diagram
$$
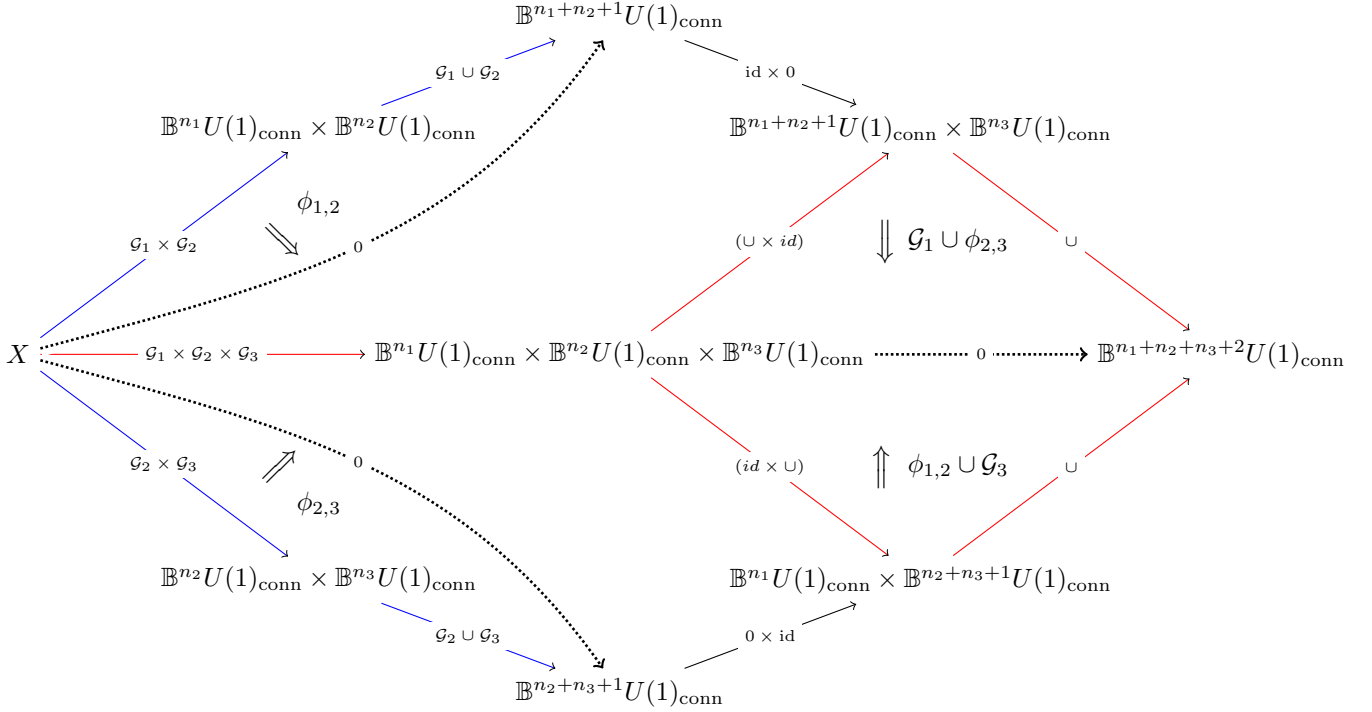
\hspace{-5mm}
\begin{tikzpicture}[descr/.style={fill=white,inner sep=2.5pt}]
\node (1) at (0,0) {\small $\BB^{n_{1}}U(1)_{\rm conn}\times \BB^{n_{2}}U(1)_{\rm conn}\times \BB^{n_{3}}U(1)_{\rm conn}$}; 
\node (2) at (4,3) {\small $\BB^{n_{1}+n_{2}+1}U(1)_{\rm conn}\times \BB^{n_{3}}U(1)_{\rm conn}$};
\node (3) at (4,-3) {\small $\BB^{n_{1}}U(1)_{\rm conn}\times \BB^{n_{2}+n_{3}+1}U(1)_{\rm conn}$};
\node (4) at (8,0) {\small $\BB^{n_{1}+n_{2}+n_{3}+2}U(1)_{\rm conn}$};
\node (5) at (-8,0) {\small $X$};
\node (6) at (-4,3) {\small $\BB^{n_{1}}U(1)_{\rm conn}\times \BB^{n_{2}}U(1)_{\rm conn}$};
\node (7) at (-4,-3) {\small $\BB^{n_{2}}U(1)_{\rm conn}\times \BB^{n_{3}}U(1)_{\rm conn}$};
\node (8) at (0,4.5) {\small $\BB^{n_{1}+n_{2}+1}U(1)_{\rm conn}$ };
\node (9) at (0,-4.5) {\small $\BB^{n_{2}+n_{3}+1}U(1)_{\rm conn}$};
\node (10) at (-4,2) {\small $\phi_{1,2}$};
\node (11) at (-4,-2) {\small $\phi_{2,3}$};
\node (12) at (-4.5,1.5) {\rotatebox{315}{$\Longrightarrow$}};
\node (13) at (-4.5,-1.5) {\rotatebox{45}{$\Longrightarrow$}};
\node (14) at (3.5,1.5) {\rotatebox{270}{$\Longrightarrow$}}; 
\node (15) at (3.5,-1.5) {\rotatebox{90}{$\Longrightarrow$}}; 
\node (16) at (4.5,1.5) {\small ${\cal G}_1\cup\phi_{2,3}$};
\node (17) at (4.5,-1.5) {\small $\phi_{1,2}\cup {\cal G}_3$};

\draw[->,draw=red] (1) to node[descr] {\tiny $(\cup\times id)$} (2);
\draw[->,draw=red] (1) to node[descr] {\tiny $(id \times \cup)$}(3);
\draw[->,draw=red] (2) to node[descr] {\tiny $\cup $}(4);
\draw[->,draw=red] (3) to node[descr] {\tiny $\cup $}(4);
\draw[->,draw=red] (5) to node[descr] {\tiny $\mathcal{G}_{1}\times \mathcal{G}_{2} \times \mathcal{G}_{3}$} (1);
\draw[->,draw=blue] (5) to node[descr] {\tiny $\mathcal{G}_{1}\times \mathcal{G}_{2}$} (6);
\draw[->,draw=blue] (5) to node[descr] {\tiny $\mathcal{G}_{2}\times \mathcal{G}_{3}$} (7);
\draw[->, draw=blue] (6) to node[descr] {\tiny $\mathcal{G}_{1}\cup \mathcal{G}_{2}$} (8);
\draw[->,draw=blue] (7) to node[descr] {\tiny $\mathcal{G}_{2}\cup \mathcal{G}_{3}$} (9);
\draw[->] (8) to node[descr] {\tiny ${\rm id}\times 0$} (2);
\draw[->] (9) to node[descr] {\tiny $0\times {\rm id}$} (3);

\draw[densely dotted,-,line width=6pt,draw=white] (5) to[out=15, in=235] node[descr] {\tiny $0$} (8);
\draw[densely dotted,->,line width=1pt] (5) to[out=15, in=235] node[descr] {\tiny $0$} (8);
\draw[densely dotted,-,line width=6pt,draw=white] (5) to[out=-15, in=125] node[descr] {\tiny $0$} (9);
\draw[densely dotted,->,line width=1pt] (5) to[out=-15, in=125] node[descr] {\tiny $0$} (9);
\draw[densely dotted,->,line width=1pt] (1) to node[descr] {\tiny $0$} (4);
\end{tikzpicture} 
$$
These two choices of homotopies $\phi_{1,2}$ and $\phi_{2,3}$ make the entire diagram homotopy commutative, as the triple cup products (the two red arrows) are trivialized by the homotopies ${\cal G}_1\cup \phi_{1,2}$ and $\phi_{2,3}\cup {\cal G}_3$. Since the cup product is strictly associative, the diagram in red commutes and we have the homotopy commuting diagram
\(
\begin{tikzpicture}[descr/.style={fill=white,inner sep=2.5pt}] 
\node (4) at (3,0) {\small $\BB^{n_{1}+n_{2}+n_{3}+2}U(1)_{\rm conn}$~.};
\node (5) at (-3,0) {\small $X$};
\node (12) at (-0.5,0.6) {\rotatebox{270}{$\Longrightarrow$}};
\node (15) at (-0.5,-0.6) {\rotatebox{90}{$\Longrightarrow$}};
\node (1) at (0.5,0.5) {\footnotesize ${\cal G}_1\cup \phi_{2,3}$}; 
\node (1) at (0.5,-0.5) {\footnotesize $\phi_{1,2}\cup {\cal G}_3$};

\draw[densely dotted,->,line width=1pt] (5) to[out=35, in=145] node[descr] {\footnotesize $0$} (4);
\draw[densely dotted,->,line width=1pt] (5) to[out=-35, in=-145] node[descr] {\footnotesize $0$} (4);
\draw[->] (5) to node[descr] {\footnotesize ${\cal G}_1\cup {\cal G}_2\cup {\cal G}_3$} (4);
\end{tikzpicture} 
\) 
These two homotopies fit together to form a loop. Then, by Prop. \ref{Prop loop} and the universal property of the homotopy pullback, 
we can equivalently describe this as a map
\(
X\longrightarrow \Omega\BB^{n_{1}+n_{2}+n_{3}+2}U(1)_{\rm conn}\simeq \flat \BB^{n_{1}+n_{2}+n_{3}+1}U(1)\;.
\label{loop}
\) 
\begin{lemma}
The homotopy class of the loop \eqref{loop} is in the image of the inclusion of the group 
$H^{n_{1}+n_{2}+n_{3}+1}(X;U(1))$ into $\widehat{H}^{n_{1}+n_{2}+n_{3}+2}(X;\ZZ)$.
\end{lemma}
\theproof
Using the Dold-Kan adjunction along with a {\v C}ech resolution of $X$, we have the following sequence of isomorphisms
\begin{eqnarray*}
\pi_{0}\map(X,\flat\BB^{n_{1}+n_{2}+n_{3}+1}U(1)) &\simeq & H_{0}\hom_{\chp}(N(C(\{U_{i}\}),{\rm disc}(U(1))[n_{1}+n_{2}+n_{3}+1])
\\
&\simeq& H^{n_{1}+n_{2}+n_{3}+1}(X;U(1))
\\
&\into& \hat{H}^{n_{1}+n_{2}+n_{3}+2}(X;\Z)\;.
\end{eqnarray*}
Here ${\rm disc}$ indicates that we are taking the discrete topology on $U(1)$ (i.e. this is the sheaf of locally constant $U(1)$ valued functions).
\endofproof
\begin{remark}
{\bf (i)} Notice that we could have equivalently taken the homotopy class of the loop directly to get an element $[\mathcal{G}_1\cup \phi_{2,3}-\phi_{1,2}\cup \mathcal{G}_1]$ in
\begin{eqnarray*}
\pi_{1}\map(X,\BB^{n_{1}+n_{2}+n_{3}+2}U(1)_{\rm conn}) &\simeq& 
H_{1}\hom_{\chp}\big(C(\{U_i\}),\ \Z^{\infty}_{\mathcal{D}}[n_{1}+n_{2}+n_{3}+2]\big)
\\
&\simeq& H^{n_{1}+n_{2}+n_{3}+1}(X;U(1))
\\
&\into& \widehat{H}^{n_{1}+n_{2}+n_{3}+2}(X;\Z)\;.
\end{eqnarray*}

\noindent {\bf (ii)} The above observations allow us to recover the usual definition of the Massey product as an element in cohomology. In section \ref{Sec Massey hyper}, we observed that such a class is not completely well-defined purely at the level of cohomology and there was some dependence on the chosen cochain representatives. Taking this point of view, one can see this dependence as a choice of trivializations $\phi_{1,2}$ and $\phi_{2,3}$ of the cup products.
\end{remark}

\medskip
This definition works well for the triple product and gives a clear picture on how the triple product is built out of the homotopies. However, to describe the higher triple products this way would be cumbersome. Moreover, the algebraic nature of the products would not be transparent. For these reasons, we will use the language of simplicial homotopy theory to describe these homotopy commuting diagrams and the Dold-Kan correspondence to organize these homotopies in an algebraic way. To prepare the reader for this perspective, we first recast the triple product in this language. 

\medskip
Notice that the triple product was described by two homotopies $\phi_{1,2}$ and $\phi_{2,3}$ connecting the basepoint $0$ to the double cup products. We can express this situation diagrammatically via
the horn-fillers
{\small 
\begin{minipage}[t]{0.35\textwidth}
\begin{equation*}
\begin{tikzpicture}
\matrix(m) [matrix of math nodes, column sep=5em, row sep=3em]{
\partial \Delta[1] & \left[X, \BB^{n_1+n_2+1}U(1)_{\rm conn}\right]\;,
\\
\Delta[1] & 
\\
};

\path[->]
(m-1-1) edge node[above] {\tiny $(0,\cG_1\cup \cG_2)$} (m-1-2)
(m-1-1) edge (m-2-1)
(m-2-1) edge[densely dotted] node[below] {\small $\phi_{1,2}$} (m-1-2);
\end{tikzpicture}
\end{equation*}
\end{minipage}
\qquad \quad
\begin{minipage}[t]{0.45\textwidth}
$$
\begin{tikzpicture}[baseline=(current bounding box.center)]
\matrix(m) [matrix of math nodes, column sep=5em, row sep=3em]{
\partial \Delta[1] & \left[X, \BB^{n_2+n_3+1}U(1)_{\rm conn}\right]\;.
\\
\Delta[1] & 
\\
};

\path[->]
(m-1-1) edge node[above] {\tiny $(0,\cG_2\cup \cG_3)$} (m-1-2)
(m-1-1) edge (m-2-1)
(m-2-1) edge[densely dotted] node[below] {\small $\phi_{2,3}$} (m-1-2);
\end{tikzpicture}
$$
\end{minipage}
}

\vspace{1mm}
\noindent Now we would like to use these homotopies to construct a loop. To do this, we need to manipulate 
algebraically  these homotopies. This motivates us to take the Moore complex of these 
diagrams in order to translate the data 
 into the language of sheaves of chain complexes. This gives the data
  
{\hspace{-9mm}  
\begin{minipage}[t]{0.45\textwidth}
$$
\begin{tikzpicture}
\matrix(m) [matrix of math nodes, column sep=5em, row sep=3em]{
\ZZ\oplus \ZZ & N\left([X,\BB^{n_1+n_2+1}U(1)_{\rm conn}]\right)_0
\\
\ZZ & N\left([X,\BB^{n_1+n_2+1}U(1)_{\rm conn}]\right)_{1}\;,
\\
};

\path[->]
(m-1-1) edge node[above] {\tiny $(0,\cG_1\cup \cG_2)$} (m-1-2)
(m-2-1) edge node[right] {\footnotesize $(1,-1)$} (m-1-1)
(m-2-2) edge node[right] {\footnotesize $\partial$} (m-1-2)
(m-2-1) edge node[above] {\footnotesize $\phi_{1,2}$} (m-2-2);
\end{tikzpicture}
$$
\end{minipage}
\qquad \quad 
\begin{minipage}[t]{0.5\textwidth}
$$
\begin{tikzpicture}[baseline=(current bounding box.center)]
\matrix(m) [matrix of math nodes, column sep=5em, row sep=3em]{
\ZZ\oplus \ZZ & N\left([X,\BB^{n_2+n_3+1}U(1)_{\rm conn}]\right)_{0}
\\
\ZZ & N\left([X,\BB^{n_2+n_3+1}U(1)_{\rm conn}]\right)_{1}\;,
\\
};

\path[->]
(m-1-1) edge node[above] {\tiny$(0,\cG_2\cup \cG_3)$} (m-1-2)
(m-2-1) edge node[right] {\footnotesize $(1,-1)$} (m-1-1)
(m-2-2) edge node[right] {\footnotesize $\partial$} (m-1-2)
(m-2-1) edge node[above] {\footnotesize $\phi_{2,3}$} (m-2-2);
\end{tikzpicture}
$$
\end{minipage}
}
\noindent where the subindices indicate the degree of the chain complex. Now we can represent these chain homotopies succinctly in the upper triangular matrix
\begin{displaymath}
A= \left( \begin{array}{cccc}
0 & \cG_1 & \phi_{1,2} & \ast \\
0 & 0 & \cG_2 & \phi_{2,3} \\
0 & 0 & 0 & \cG_3 \\
0 & 0 & 0 & 0
\end{array} \right)\;.
\end{displaymath}
By construction, this matrix satisfies the Maurer-Cartan equation
$$
dA-\overline{A}\cdot A=\mu(A)\in \rm{Ker(A)}\;.
$$
Moreover, $\mu(A)$ is of the form
\begin{displaymath}
\mu(A)= \left( \begin{array}{cccc}
0 & 0 & 0 & \tau \\
0 & 0 & 0 & 0 \\
0 & 0 & 0 & 0 \\
0 & 0 & 0 & 0
\end{array} \right)\;.
\end{displaymath}
Applying the differential $d$ to $\tau$ and using the Leibniz rule, we get
\begin{eqnarray*}
d(\tau) &=& d\left(\cG_1\cup \phi_{2,3}-\phi_{1,2}\cup \cG_3\right)
\\
&=&  d(\cG_1)\cup \phi_{2,3}+\cG_1\cup d(\phi_{2,3})-d(\phi_{1,2})\cup \cG_3+\phi_{1,2}\cup d(\cG_3)
\\
&=& \cG_1\cup (\cG_2\cup \cG_3)-(\cG_1\cup \cG_2)\cup \cG_3
\\
&=& 0\;.
\end{eqnarray*}
At the level of sheaf hypercohomology, we have the following:
\begin{proposition} The cohomology class of the matrix cocycle $\mu (A)$ is the element  
$$
[\mu(A)] = \left[\cG_1\cup \phi_{2,3}-\phi_{1,2}\cup \cG_3\right]
\in H^{n_{1}+n_{2}+n_{3}+1}(X;U(1))\;.
$$
\end{proposition}
\theproof
We have the following sequence of isomorphisms
\begin{eqnarray*}
[\mu(A)] = \left[\cG_1\cup \phi_{2,3}-\phi_{1,2}\cup \cG_3\right]
&\in& H_{1}\hom_{\chp}\left(C(\mathcal{U}),\ \BB^{n_{1}+n_{2}+n_{3}+2}U(1)_{\rm conn}\right)
\\
&&\simeq \pi_{1}\map\left(C(\mathcal{U}),\ \BB^{n_{1}+n_{2}+n_{3}+2}U(1)_{\rm conn}\right)
\\
&&\simeq \pi_{0}\map\left(C(\mathcal{U}),\ \flat\BB^{n_{1}+n_{2}+n_{3}+1}U(1)\right)
\\
&&\simeq  H^{n_{1}+n_{2}+n_{3}+1}(X;U(1))\;.
\end{eqnarray*}

\vspace{-5mm}
\endofproof

\subsection{General stacky Massey products}

\label{Sec General stacky} 


We would like to utilize the machinery of May \cite{May} which makes use of matrices. 
We will introduce stacks labelled by two integers, which will be indexing the entries
of the corresponding matrices. 
To that end, let $\mathcal{R}_{ij}$, $i,j \in \N$, be simplicial abelian stacks equipped with maps
$$
\cup:\mathcal{R}_{ij}\otimes \mathcal{R}_{jk}\longrightarrow \mathcal{R}_{ik}\;,
$$
which are associative in the sense that
$\cup \circ (\cup\otimes \mathrm{id})=\cup \circ (\mathrm{id}\otimes \cup)$. 

\begin{remark}
Let $N$ denote the normalized Moore functor.  It follows from the definition of the 
differential on the tensor product that the induced product
$$
\tilde{\cup}:N(\mathcal{R}_{ij})\otimes N(\mathcal{R}_{jk})\overset{\sim}{\longrightarrow}
 N(\mathcal{R}_{ij}\otimes \mathcal{R}_{jk})
\longrightarrow N(\mathcal{R}_{ik})
$$
must satisfy the Leibniz type rule
$
d(\alpha\cup \beta)=d(\alpha)\cup \beta+(-1)^{\mathrm{deg}}\alpha\cup d(\beta)
$
on sections.
\end{remark}
We can now utilize an extension of the machinery of May \cite{May} locally to define 
the refined matric Massey products in our setting. 
To this end, we consider the set of all upper triangular half-infinite matrices
$M(\cR)=\bigcup_n M(\cR)_n$, where (cf. \eqref{MA}) 
\(\label{matrices}
M(\cR)_n=\left\{ A=(a_{ij}) ~\vert~ a_{ij}\in N(\cR_{ij}), a_{ij}=0 ~{\rm for}~
j\leq i~ {\rm and}~ i,j \geq n+1 ~{\rm for~ some}~ n \in \N \right\}
\)
is the subalgebra of $n \times n$ matrices.  
Notice that, with our definition, this set possesses more structure. It becomes a sheaf of 
DGA's with product given by matrix multiplication and differential given by applying 
the differential on $N(\cR_{ij})$ to each entry of the matrix. Just as in the case of 
classical Massey products, we have a filtration of presheaves of subalgebras
\(
M(\cR)_{1}\subset M(\cR)_{2}\subset \hdots \subset M(\cR)_{n}\subset \hdots 
\)
and a bigrading
\(
M(\cR)=\sum _{p\geq 1, k\geq 0}M^{p,k}\;,
\)
where 
\(
M^{p,k}={\rm span}\left \{ 
{\footnotesize
\left(\begin{array}{ccccc} 
& & & & 0
\\
& & a_{i,i+p} & &
\\
& & & & 
\\
& 0 & & &
\\
& & & &
\end{array}\right) 
}
; a_{i,i+p}\in N(\cR_{i,i+p})\right\}.
\)
\noindent We can define the following notions similarly to the classical case.  
\begin{definition}\label{maurer1}
Let $A$ be a matrix in $M(\cR)$. We define the (stacky version) of the 
Maurer-Cartan equation as
$$
dA-\overline{A}\cdot A\equiv 0~~{\rm mod} \ker(A)\;,
$$
and call a solution a \emph{formal} connection with \emph{curvature}
$$
\mu(A)=dA-\overline{A}\cdot A\;.
$$
\end{definition}

We are now ready to define the stacky Massey product with a product on the bigraded 
sequence of stacks. 

\begin{definition} 
\label{Massey2} 
Let $\cR=\{\cR_{ij}\}$ be a sequence of abelian stacks equipped with maps
$$
\cup:\cR_{ij}\otimes \cR_{jk}\longrightarrow \cR_{ik}\;,
$$
which satisfy
$
\cup\circ (\rm{id}\otimes \cup)=\cup \circ (\cup \otimes \rm{id})
$.
Let $A$  be a formal connection with curvature $\mu(A)$. Then the entries of the hypercohomology class $[\mu(A)]$ are called \emph{stacky Massey products}.
\end{definition} 

\begin{remark}
The following examples of stacks satisfy the compatibility requirement of Def. \ref{Massey2} and will be of particular interest to us.
They are the mapping stacks corresponding to the stacks described in Remark \ref{stacks}.
Fix a manifold $X$ and a sequence $(n_{i,j})$, $i<j\leq n$, of integers satisfying $n_{i,j}+n_{j,k}=n_{i,k}$;
\begin{enumerate} [label= {\bf (\roman*)}]
\item The stacks $[X,\BB^{n_{i,j}-1}U(1)_{\rm conn}]$ of higher bundles with connection, with the stacky cup product and
 {\v C}ech-Deligne differential.
\item The stacks $[X,\BB^{n_{i,j}}\ZZ]$ of higher bundles, with the usual cup product and singular differential.
\item The stacks $[X,\flat_{\rm dR}\BB^{n_{i,j}}U(1)_{\rm conn}]$ of differential forms of degrees $\leq n$, with the wedge 
product and exterior derivative.
\end{enumerate}
\end{remark}

We highlight the power of the above definitions in the following examples, where we  are able to 
describe all three of the differential, singular, and de Rham triple products.

\begin{example} {\bf (Differential triple product)} Let $\cG_{i}$, $i=1,2,3$, be bundles 
corresponding to morphisms $\Delta[0]\to [X,\BB^{n_{i,i+1}-1}U(1)_{\rm conn}]$. Suppose 
$\cG_{1}\cup \cG_{2}$ and $\cG_{2}\cup \cG_{3}$ represent trivial classes in 
$\pi_{0}\map(X,\BB^{n_{i,j}-1}U(1)_{\rm conn})$. 
Choose a defining system
\begin{displaymath}
A= \left( \begin{array}{cccc}
0 & \cG_1 & \phi_{1,2} & \ast \\
0 & 0 & \cG_2 & \phi_{2,3} \\
0 & 0 & 0 & \cG_3 \\
0 & 0 & 0 & 0
\end{array} \right)\;,
\end{displaymath}
where $\phi_{1,2}$ and $\phi_{2,3}$ are nondegenerate 1-simplices trivializing the cup products. Then 
the curvature of $A$ is 
$$
\mu(A)= \left( \begin{array}{cccc}
0 & 0 & 0 & \tau \\
0 & 0 & 0 & 0 \\
0 & 0 & 0 & 0 \\
0 & 0 & 0 & 0
\end{array} \right)\;,
$$
and the hypercohomology class $[\tau]$ is 
$
\left[\cG_1\cup \phi_{2,3}-\phi_{1,2}\cup \cG_3\right]
$.
The latter is an element in
$$
\pi_{1}\map(X,\BB^{n_{3,3}-1}U(1)_{\rm conn}) \simeq  \pi_{0}\map(X,\flat\BB^{n_{1,4}-2}U(1)) \simeq H^{n_{1,4}-2}(X;U(1))\;,
$$
where we have $n_{1,4}=n_{1,3}+n_{3,4}=n_{1,2}+n_{2,3}+n_{3,4}$.
\end{example}

\begin{example} {\bf(Singular triple product)} Let $X$ be a manifold, and let $\vert X \vert$ be the topological 
space denoting its geometric realization. Let $a_{i}:\vert X \vert \to K(\ZZ, n_{i,i+1})\simeq B^{n_{i,i+1}}\ZZ$, 
$i=1,2,3$, be singular cochains with cup products vanishing in cohomology. Choose a defining system
\begin{displaymath}
A= \left( \begin{array}{cccc}
0 & a_1 & f_{1,2} & \ast \\
0 & 0 & a_2 & f_{2,3} \\
0 & 0 & 0 & a_3 \\
0 & 0 & 0 & 0
\end{array} \right)\;.
\end{displaymath}
Since geometric realization is a left $\infty$-adjoint the discrete stack functor $\mathrm{disc}$ \cite{Urs}, 
these are equivalently given by maps of stacks
$$
\bar{a}_{i}:\Delta[0]\longrightarrow [X, \BB^{n_{i,i+1}}\ZZ]\;,
$$
and homotopies
$$
\bar{f}_{i,i+1}:\Delta[1] \longrightarrow [X, \BB^{n_{i,i+2}}\ZZ]
$$
trivializing the cup products, hence a defining system
\begin{displaymath}
A= \left( \begin{array}{cccc}
0 & \bar{a}_1 & \bar{f}_{1,2} & \ast \\
0 & 0 & \bar{a}_2 & \bar{f}_{2,3} \\
0 & 0 & 0 & \bar{a}_3 \\
0 & 0 & 0 & 0
\end{array} \right)\;.
\end{displaymath} 
The hypercohomology class of the entry $\tau\in \mu(A)$ is given by
$
\left[\bar{a}_1\cup \bar{f}_{2,3}-\bar{f}_{1,2}\cup \bar{a}_3\right]
$,
which is an element in
\bea
\pi_{1}\map(X,\BB^{n_{1,4}}\ZZ)&\simeq & \pi_{1}\map(\vert X \vert, K(\ZZ,n_{1,4}))
\nonumber\\
&\simeq & \pi_{0}\map(\vert X \vert, K(\ZZ,n_{1,4}-1))
\nonumber\\
&\simeq& H^{n_{1,4}-1}(X,\ZZ)\;.
\eea
\end{example}

\begin{example} {\bf(de Rham triple product)}
Let $X$ be a manifold and let $\alpha_{i}$, $i=1,2,3$, be closed forms in different degrees. 
These forms are equivalently given by maps
$$
\alpha_{i}:\Delta[0]\longrightarrow \flat_{\rm dR}\BB^{n_{i,i+1}}U(1)_{\rm conn}\;.
$$
Suppose that the wedge products 
$\alpha_{1}\wedge \alpha_2$ and $\alpha_2\wedge \alpha_3$ are trivial in cohomology. 
Then we can choose a defining system via
\begin{displaymath}
A= \left( \begin{array}{cccc}
0 & \alpha_1 & \eta_{1,2} & \ast \\
0 & 0 & \alpha_2 & \eta_{2,3} \\
0 & 0 & 0 & \alpha_3 \\
0 & 0 & 0 & 0
\end{array} \right)\;,
\end{displaymath}  
where $\eta_{1,2}$ and $\eta_{2,3}$ are 1-simplices. The hypercohomology class of the 
entry $\tau\in \mu(A)$ is given by
$$
\left[\alpha_1\wedge \eta_{2,3}-\eta_{1,2}\cup \alpha_3\right]\;.
$$
The sheaf at each level in the complex $\Omega^{\leq n_{1,4}}$ is acyclic (the sheaves are that of differential forms and so admit a partition of unity). Thus, we can calculate the hypercohomology as
\bea
\pi_{1}\map(X, \Omega^{\leq n_{1,4}})&\simeq &
  H_{1}\Omega^{\leq n_{1,4}}(X)
  \nonumber\\
  &\simeq &H_{\rm dR}^{n_{1,4}-1}(X)\;.
\eea
\end{example}

Our main result in this section relates Massey products for Deligne cocycles to 
corresponding ones for higher bundles in the stacky sense. 

\begin{theorem}
Let $\hat{a}_{i}$, $1\leq i \leq l$, be Deligne cocycles. Suppose the $l$-fold Massey product is defined. 
Let $\cG_i$, $1\leq i\leq l$, be $n_{i, i+1}$- bundles with connections
$$
\cG_{i}:X\longrightarrow \BB^{n_{i,i+1}}U(1)_{\rm conn}\;,
$$
representing the Deligne cocycles. Then there is a natural bijection between corresponding Massey products
$$
\langle \cG_{1},\cG_{2},\hdots ,\cG_{l} \rangle \simeq \langle \hat{a}_{1},\hat{a}_{2},\hdots ,\hat{a}_{l} \rangle\;.
$$
\end{theorem}
\theproof
Recall that $\BB^{n}U(1)_{\rm conn}:=\Gamma(\ZZ^{\infty}_{\mathcal{D}}[n+1])$ (see \cite{FSSt}).
 Using the definition of the
 stacky hom, the fact that the counit $\epsilon:N\Gamma\to {\rm id}$ is a natural isomorphism and the lax monoidal structure on $N$, we have a homotopy equivalence for each test object $U$, 
\begin{eqnarray*}
N([X,\BB^{n}U(1)_{\rm conn}])(U) &=& N([X,\BB^{n}U(1)_{\rm conn}](U))
\\
&\simeq & N(\map(C(\{U_i\})\times U,\BB^{n}U(1)_{\rm conn}))
\\
&\simeq & \hom_{\chp}(N(C(\{U_i\})\otimes N(U),\ZZ^{\infty}_{\mathcal{D}}[n+1])
\\
&\simeq&  \hom_{\chp}(N(C(\{U_i\}),\ZZ^{\infty}_{\mathcal{D}}[n+1](U))
\\
&\simeq& C(X,\ZZ^{\infty}_{\mathcal{D}}[n+1])(U)\;,
\end{eqnarray*}
where the last line denotes the {\v C}ech resolution of the Deligne complex $\ZZ^{\infty}_{\mathcal{D}}[n+1]$. Hence, a defining system in the stacky sense is naturally equivalent to a defining system in the sense of \cite{Sch}. 
Since the set of Massey products is parametrized by the set of defining systems, it follows that indeed we have a natural bijection
$$
\langle \cG_{1},\cG_{2},\hdots ,\cG_{l} \rangle \simeq \langle \hat{a}_{1},\hat{a}_{2},\hdots ,\hat{a}_{l} \rangle\;.
$$

\vspace{-5mm}
\endofproof

\subsection{Properties of stacky Massey products} 

\label{Sec properties} 


We will now consider properties of the stacky Massey products. Our setting allows for these to 
be quite attractive and natural. The most immediate of those are direct generalizations  of classical 
ones. Later in this section we will see properties that are more peculiar to the differential 
setting. Among the properties that the
 classical Massey products satisfy are the following (see \cite{May} \cite{Kr}):
\begin{enumerate}[label= {\bf (\roman*)}]
\item {\it Dimension:} The dimension of $\langle x_1, x_2, \cdots, x_n \rangle$ is $\sum \deg (x_i) - n +2$.
\item {\it Naturality:} If $f: X \to Y$ is a continuous map and $y_1 \cdots, y_k \in H^*(Y; R)$ 
such that the $k$-fold Massey product $\langle y_1, y_2, \cdots, y_k \rangle$ is defined, 
then $\langle x_1, \cdots, x_k \rangle= \langle f^*(y_1), \cdots, f^*(y_k) \rangle$ is defined 
as a Massey product on the cohomology of $X$ and 
$$
f^*(\langle y_1, \cdots, y_k \rangle) \subset \langle f^*(y_1), \cdots, f^*(y_k) \rangle\;.
$$
\item {\it Definedness:} The vanishing of the the lower Massey products is only a necessary condition for the $k$-fold Massey product to be defined for $k>3$. For $k=3$ the condition is both necessary and sufficient.

\item {\it Slide relation:} If the Massey product $\langle x_1, x_2,\hdots, x_n \rangle$ is defined, then so is $\langle x_1, x_2,\hdots, rx_i, \hdots x_n \rangle$ for any $r\in R$. Moreover we have the relation 
$$r\langle x_1, x_2,\hdots , x_n \rangle\subset \langle x_1,x_2,\hdots, rx_i, \hdots x_n \rangle\;.$$ 
\end{enumerate}

These indeed extend to the stacky version.

\begin{proposition}
The stacky Massey products satisfy the following properties:
\begin{enumerate}[label= {\bf (\roman*)}]
\item {\rm Dimension:}  The dimension of $\langle \mathcal{G}_1, \mathcal{G}_2, \cdots, \mathcal{G}_l \rangle$ is $\sum \deg (\mathcal{G}_i) - l +2$.

\item  {\rm Naturality}: If $f: X \to Y$ is a smooth map between manifolds and $\mathcal{G}_1 \cdots, \mathcal{G}_k \in \hat{H}_{\mathcal{D}}^*(X ; \ZZ)$ 
such that the $k$-fold Massey product $\langle \mathcal{G}_1, \mathcal{G}_2, \cdots, \mathcal{G}_k \rangle$ is defined, 
then $\langle \mathcal{G}_1, \cdots, \mathcal{G}_k \rangle= \langle f^*(\mathcal{G}_1), \cdots, f^*(\mathcal{G}_k) \rangle$ is defined 
as a Massey product on the differential cohomology of $X$ and 
$$
f^*(\langle \mathcal{G}_1, \cdots, \mathcal{G}_k \rangle) \subset \langle f^*(\mathcal{G}_1), \cdots, f^*(\mathcal{G}_k) \rangle\;.
$$
\item {\rm Definedness}: The vanishing of the the lower Massey products is only a necessary condition for the $k$-fold Massey product to be defined for $k>3$. For $k=3$ the condition is both necessary and sufficient.

\item {\rm Slide relation:} If the Massey product $\langle {\cal G}_1, {\cal G}_2,\hdots, {\cal G}_n \rangle$ is defined, then so is $\langle {\cal G}_1, {\cal G}_2,\hdots, m{\cal G}_i, \hdots {\cal G}_n \rangle$ for any $m\in \ZZ$. Moreover we have the relation 
$$m\langle {\cal G}_1, {\cal G}_2,\hdots , {\cal G}_n \rangle\subset \langle {\cal G}_1,{\cal G}_2,\hdots, m{\cal G}_i, \hdots {\cal G}_n \rangle\;.$$ 
\end{enumerate}

\end{proposition}

 \theproof
Part 1 follows immediately from the definition. To prove part 2, note that the functor $[-,\cR]$ is contravariant, sending a map $f:X\to Y$ to its pullback
$$
f^{*}:[Y,\cR_{ij}]\longrightarrow [X,\cR_{ij}]\;.
$$
Since the cup product is natural with respect to pullbacks, the induced morphism $f^*:N([Y,\cR_{ij}]) \to N([X,\cR_{ij}])$
descends to a morphism of sheaves of DGA's
$$
f^*:M([Y,\cR_{ij}])\longrightarrow M([X,\cR_{ij}])\;.
$$
It follows that if $A$ is a formal connection in $M([Y,\cR_{ij}])$, then $f^*(A)$ is a formal connection in $M([X,\cR_{ij}])$ satisfying the equation:
$$
df^*(A)-f^*(A)\cdot \overline{f^*(A)}=f^*(\mu(A))\in \ker(f^*(A))\;.
$$
By definition of the $k$-fold Massey product, the claim follows. For part 3, we will show that for $k=3$ the condition is both necessary and sufficient. From the proof, it will be clear that this cannot be the case for higher products. Let $\cG_1$, $\cG_2$ and $\cG_3$ be bundles and suppose the triple product $\langle \cG_1,\cG_2,\cG_3\rangle$ is defined. Then we have trivializations $\phi_{1,2}$ and $\phi_{2,3}$ such that
$$
d\phi_{1,2}=\cG_1\cup \cG_2 \qquad \text{and} \qquad d\phi_{2,3}=\cG_2\cup \cG_3\;.
$$
Hence, both cup products are trivial. For the converse, it is clear that if both cup products are trivial in cohomology, we can choose trivializing homotopies and form the Massey triple product. For higher products, the higher trivializations depend on the lower ones. In fact, for the fourfold product, choose trivializations $\phi_{1,2}$, $\phi_{2,3}$ and $\phi_{3,4}$ of the cup products such that
$$\cG_1\cup \phi_{2,3}-\phi_{1,2}\cup \cG_3\in \langle \cG_1,\cG_2,\cG_3\rangle$$
is trivializable. Then for the fourfold Massey product to be defined, the other triple product
$$\cG_2\cup \phi_{2,3}-\phi_{3,4}\cup \cG_4\in \langle \cG_2,\cG_3,\cG_4\rangle$$
must be trivializable. But this may not be true, even if $\langle \cG_2,\cG_3,\cG_4\rangle$ contains $0$. Finally, for part 4, let $A$ be a formal connection of the form
$$
{\small 
A=
\left( \begin{array}{cccccccccc}
 0 & \cG_1 & \phi_{1,2} & \hdots &&&&&& \ast
\\
 &     0   & \cG_2 & \phi_{2,3} & \hdots &&&&&
\\
  &        &     0    & \hdots & \hdots & \hdots  &&&
  \\
  &        &          &             & \cG_{i-1} & \phi_{i-1,i}&&&&
  \\
  &         &          &            &             & \cG_i & \phi_{i, i+1} &&&
  \\
  &         &           &            &           &             & \cG_{i+1} &&&
  \\
  &         &          &            &              &           &                & \hdots &\phi_{n-2,n-1}&
\\
  &        &          &             &              &           &                 &            & \cG_{n-1} & \phi_{n-1,n}
 \\ 
  &         &          &             &           &            &                 &             &    & \cG_{n}
  \\
  &&&&&&&&& 0
  \end{array}
  \right)
  }
  \;.
  $$
Then the matrix
$$
{\small 
\tilde{A}=
\left( \begin{array}{ccccc|ccccc}
 0 & \cG_1 & \phi_{1,2} & \hdots & \hdots & m\phi_{1,i} & \hdots &&& \ast
\\
 &     0   & \cG_2 & \phi_{2,3} & \hdots & \hdots &&&&
\\
  &        &     0    & \hdots & \hdots & m\phi_{i-2,i} & \hdots &&
  \\
  &        &          &             & \cG_{i-1} & m \phi_{i-1,i}& m\phi_{i-1,i+1}& \hdots &&
  \\
  &         &          &            &             & m\cG_i & m\phi_{i, i+1} & m\phi_{i,i+2}& \hdots & m\phi_{i, n} 
  \\ \hline
  &         &           &            &           &             & \cG_{i+1} & \hdots &
  \\
  &         &          &            &              &           &                & \hdots &\phi_{n-2,n-1}
\\
  &        &          &             &              &           &                 &         0  & \cG_{n-1} &  \phi_{n-1,n}
 \\ 
  &         &          &             &           &            &                 &             &   0             & \cG_{n}
  \\
  &&&&&&&&& 0
  \end{array}
  \right)
  }
  $$
is also a formal connection: that is, a defining system for the Massey product $\langle {\cal G}_1, \hdots, m{\cal G}_i, \hdots, {\cal G}_n \rangle$. Indeed, let us write the matrix $A$ as a block matrix
$$
{\footnotesize
A=
\left( \begin{array}{ccc|ccc}
 & \raisebox{-8pt}{\mbox{$A_1$}} &&&   \raisebox{-8pt}{\mbox{$A_2$}}&
 \\[3ex] \hline
   &\raisebox{-8pt}{\mbox{$0$}} &&&  \raisebox{-8pt}{\mbox{$A_3$}}&
   \\[3ex]
  \end{array}
  \right)
  }
  \;.
  $$
  Then the second matrix can be written
  $$
 { \footnotesize
  \tilde{A}=
\left( \begin{array}{ccc|ccc}
 & \raisebox{-8pt}{\mbox{$A_1$}} &&&   \raisebox{-8pt}{\mbox{$mA_2$}}&
 \\[3ex] \hline
   &\raisebox{-8pt}{\mbox{$0$}} &&&  \raisebox{-8pt}{\mbox{$A_3$}}&
   \\[3ex]
  \end{array}
  \right)
  }
  \;.
  $$
  Now the Maurer-Cartan equation for $\tilde{A}$ reads
  \bea
  \mu(\tilde{A}) &=&
  {\footnotesize
\left( \begin{array}{ccc|ccc}
 & \raisebox{-8pt}{\mbox{$dA_1$}} &&&   \raisebox{-8pt}{\mbox{$mdA_2$}}&
 \\[3ex] \hline
   &\raisebox{-8pt}{\mbox{$0$}} &&&  \raisebox{-8pt}{\mbox{$dA_3$}}&
   \\[3ex]
  \end{array}
  \right) - 
\left( \begin{array}{ccc|ccc}
 & \raisebox{-8pt}{\mbox{$\overline{A_1}$}} &&&   \raisebox{-8pt}{\mbox{$\overline{mA_2}$}}&
 \\[3ex] \hline
   &\raisebox{-8pt}{\mbox{$0$}} &&&  \raisebox{-8pt}{\mbox{$\overline{A_3}$}}&
   \\[3ex]
  \end{array}
  \right)
\left( \begin{array}{ccc|ccc}
 & \raisebox{-8pt}{\mbox{$A_1$}} &&&   \raisebox{-8pt}{\mbox{$mA_2$}}&
 \\[3ex] \hline
   &\raisebox{-8pt}{\mbox{$0$}} &&&  \raisebox{-8pt}{\mbox{$A_3$}}&
   \\[3ex]
  \end{array}
  \right)
  }
  \\
  &=&
  {\footnotesize
  \left( \begin{array}{ccc|ccc}
 & \raisebox{-8pt}{\mbox{$dA_1$}} &&&   \raisebox{-8pt}{\mbox{$mdA_2$}}&
 \\[3ex] \hline
   &\raisebox{-8pt}{\mbox{$0$}} &&&  \raisebox{-8pt}{\mbox{$dA_3$}}&
   \\[3ex]
  \end{array}
  \right)
  -
  \left( \begin{array}{ccc|ccc}
 & \raisebox{-8pt}{\mbox{$\overline{A_1}A_1$}} &&&   \raisebox{-8pt}{\mbox{$m(\overline{A_1}A_2+\overline{A_2}A_3)$}}&
 \\[3ex] \hline
   &\raisebox{-8pt}{\mbox{$0$}} &&&  \raisebox{-8pt}{\mbox{$ \overline{A_3}A_3$}}&
   \\[3ex]
  \end{array}
  \right)
  }
  \;.
  \eea
   We would like to show that $\mu(\tilde{A})$ is in $\ker(\tilde{A})$. Since $A$ satisfies the Maurer-Cartan 
   equation up to an element in the kernel 
   \(
   \ker(A)=\left(\begin{array}{cccc}
   &\hdots & 0& \ast
   \\
   & \hdots & 0& 0
  
  \\
  &&\vdots& \vdots  
 \\
 &&&
 \end{array}\right)
 \)
 we must have $dA_1=\overline{A_1}\cdot A_1$ and $dA_3=\overline{A_3}\cdot A_3$. Since $A$ is a formal connection, 
 we must also have
 $$
 \mu(A)_{2}=dA_2-(\overline{A_1}A_2+\overline{A_2}A_3)
 $$
 where $\mu(A)_2$ is the upper right block of $\mu(A)$ of dimension ${\rm dim}(A_2)$. Since the only nonzero term of $\mu(A)$ is the cochain representative of the Massey product $\tau$, located in the upper right corner of $\mu(A)$, we have that
 $$
 \mu(\tilde{A})_2=mdA_2-m(A_1\overline{A_2}+A_2\overline{A_3})=m\mu(A)_2
 $$
 has one nonzero element $\sigma=m\tau$ in the upper right corner. Therefore, $\tilde{A}$ is indeed a formal connection
 and, at the level of cohomology, the only nonzero term of the class $[\mu(A)]$ is $[\sigma]=m[\tau]$. Since $[\tau]$ was 
 chosen to be an arbitrary element of the Massey product $\langle \cG_1,\hdots, \cG_n\rangle$, we have
 $$
 m\langle \cG_1,\hdots, \cG_n\rangle \subset \langle \cG_1,\hdots, m\cG_i, \hdots \cG_n\rangle\;.
 $$ 

\vspace{-6mm}
\endofproof

We now discuss the relationship between the stacky Massey product and the singular Massey product.
The following parametrizes how forgetting the differential data on the Massey product is not quite
the same as taking the Massey product of cohomology classes after forgetting the differential data on these.

\begin{proposition}
Let $\cG_{i}:\Delta[0] \to [X,\BB^{n_{i,i+1}}U(1)_{\rm conn}]$, $1\leq i\leq l$, be
higher  bundles on $X$ with defined Massey
 product. Then precomposition with the forgetful morphism
$$
I:\BB^{n}U(1)_{\rm conn}\longrightarrow \BB^{n+1}\ZZ\;,
$$
induced by the map
$$
\xymatrix{
\underline{\ZZ} \ar[d]^-0 \ar[rr]^i && \Omega^0 \ar[d]^-0 \ar[rr] && \cdots \ar[d] \ar[rr] && \Omega^{n-1} \ar[d]^-d \\
\underline{\ZZ} \ar[rr]^-d && 0 \ar[rr] && \cdots \ar[rr] && 0
\;,
}
$$
yields singular cocycles with defined Massey product. Furthermore, we have
$$
I\langle \cG_{1},\cG_{2},\cG_{3} \rangle \subset \langle I(\cG_{1}),I(\cG_{2}),I(\cG_3) \rangle\;.
$$
\end{proposition}

\theproof
For simplicitiy, we denote the sheaf of matrix algebras
\bea \label{def}
M_{\rm diff} &:=& M([X,\BB^*U(1)_{\rm conn}])
\\
M_{\rm sing} &:=& M([X,\BB^*\Z])
\eea
according to the corresponding cohomology theories for these matrices. It is clear by definition that $I$ respects the cup product structure, hence $I$ induces a morphism of sheaves of DGA's
$I_{*}:M_{\rm diff}\to M_{\rm sing}$.
It follows immediately from the definition of the Maurer-Cartan equation Def. \ref{maurer1}, 
that formal connections are sent to formal connections. Then passing to 
hypercohomology gives the result.
\endofproof

\begin{remark}
{\bf (i)} It follows from the proposition that if the classical Massey product 
$\langle I(\cG_{1}),I(\cG_{2}),I(\cG_3) \rangle$ is zero then certainly the 
left hand side is zero, i.e. $\langle \cG_{1},\cG_{2},\cG_{3} \rangle$ is in the 
kernel of the forgetful morphism $I$.
 From the sequence $\Omega^{n-1}/{\rm Im}(d) \to \hat{H}^n 
\buildrel{I}\over{\longrightarrow} H^n$ we have that $\langle \cG_{1},\cG_{2},\cG_{3} \rangle$
will be an $(n-1)$-form. 
However, it is important to note that this is {\rm not} quite 
the $(n-1)$-form given by the classical Massey product. 
\item {\bf (ii)} A related question is to ask whether the differential Massey product 
completely refines the singular Massey product. That is: 
do we have a \emph{bijection},
$$
I\langle \cG_{1},\cG_{2},\cG_{3} \rangle \simeq \langle I(\cG_{1}),I(\cG_{2}),I(\cG_3) \rangle\;?
$$
Unfortunately, this {\it cannot} be possible. Essentially, this is because the map
$I_{*}:M_{\rm diff}\to M_{\rm sing}$
has a nontrivial kernel. Hence we cannot expect the Maurer-Cartan equation to hold after refining. 
\item {\bf (iii)} However, this does help explain the nature of differential Massey products. In fact, since these products are always flat,
 it follows from diagram \eqref{diamond} that \emph{if} the refinement of a singular formal connection is again a formal connection, then the singular Massey product must have been torsion. 
 \end{remark} 
 
 We will show that the failure of the refinement to satisfy the Maurer-Cartan equation can be measured by the de Rham Massey product. 

\begin{lemma}
Let $\cF_{ij}\to \cR_{ij}\onto \mathcal{S}_{ij}$ be a fibration sequence of abelian prestacks for each $i$ and $j$. Suppose,
moreover, that we have commuting diagrams
$$
{\small 
\begin{tikzpicture}
\matrix (m) [matrix of math nodes, column sep=3em, row sep=3em]{
\cF_{ij}\otimes \cF_{jk} & \cF_{ik}
\\
\cR_{ij}\otimes \cR_{jk} & \cR_{ik}
\\
\mathcal{S}_{ij}\otimes \mathcal{S}_{jk} & \mathcal{S}_{ik}\;.
\\
};

\path[->]
(m-1-1) edge node[above] {\footnotesize $\cup$} (m-1-2)
(m-2-1) edge node[above] {\footnotesize $\cup$} (m-2-2)
(m-3-1) edge node[above] {\footnotesize $\cup$} (m-3-2)
(m-1-1) edge node[right] {\footnotesize $i\otimes i$} (m-2-1)
(m-2-1) edge node[right] {\footnotesize $p\otimes p$} (m-3-1)
(m-1-2) edge node[right] {\footnotesize $i$} (m-2-2)
(m-2-2) edge node[right] {\footnotesize $p$} (m-3-2);
\end{tikzpicture}
}
$$
Then the induced sequence $0\to M(\cF)\to M(\cR) \onto M(\mathcal{S})\to 0$ is a short exact sequence of DGA's.
\end{lemma}
\theproof
Since the normalized Moore functor is right Quillen and preserves equivalences, it follows that it sends fiber sequences to fiber sequences. Hence, we have a diagram
$$
{\small 
\begin{tikzpicture}
\matrix (m) [matrix of math nodes, column sep=3em, row sep=3em]{
N(\cF_{ij})\otimes N(\cF_{jk}) & N(\cF_{ik})
\\
N(\cR_{ij})\otimes N(\cR_{jk}) & N(\cR_{ik})
\\
N(\mathcal{S}_{ij})\otimes N(\mathcal{S}_{jk}) & N(\mathcal{S}_{ik})\;,
\\
};

\path[->]
(m-1-1) edge node[above] {\footnotesize $\cup$} (m-1-2)
(m-2-1) edge node[above] {\footnotesize $\cup$} (m-2-2)
(m-3-1) edge node[above] {\footnotesize $\cup$} (m-3-2)
(m-1-1) edge node[right] {\footnotesize $i\otimes i$} (m-2-1)
(m-2-1) edge node[right] {\footnotesize $p\otimes p$} (m-3-1)
(m-1-2) edge node[right] {\footnotesize $i$} (m-2-2)
(m-2-2) edge node[right] {\footnotesize $p$} (m-3-2);
\end{tikzpicture}
}
$$
where the right hand side is a short exact sequence of presheaves of chain complexes. By definition, it follows that we have a short exact sequence 
$$
0\to M(\cF)\to M(\cR)\onto M(\mathcal{S})\to 0
$$
of chain complexes. By commutivity of the above diagram, both maps are homomorphisms of presheaves of DGA's.
\endofproof

It follows from the lemma along with diagram \eqref{diffhex}, that there is a short exact sequence of presheaves of bigraded rings
\(\label{ses}
0 \longrightarrow M([X,\Omega^{\leq *}]) \longrightarrow M_{\rm diff}\longrightarrow M_{\rm sing}\longrightarrow 0\;.
\)
Hence, $M_{\rm form}:=M([X,\Omega^{\leq *}])$ is a two-sided ideal in $M_{\rm diff}$.

\medskip
Now, by definition of $\ker(A)$ along with the above observation, we have
$$
\ker(\hat{A})\subset \widehat{\ker(A)}\;,
$$
where $\ \widehat{\ }\ $ denotes a choice of differential refinement. In fact, for a matrix $C\in \ker(\hat{A})$ and $C^{\prime}\in \widehat{\ker(A^{\prime})}$, we have that the difference $C-C^{\prime}=B\in M_{\rm form}$. 
It is this lack of commutativity between taking kernels and taking differential refinements that
leads to a nontrivial structure than might otherwise be anticipated.

\medskip
Summarizing the previous observations gives the following theorem.
\begin{theorem}
\label{Thm MC}
Let $A$ be a formal connection for $M_{\rm sing}$, and let $\hat{A}$ be a differential refinement of $A$ with $\mu(A)$ a solution to the Maurer-Cartan equation. Then 
any differential refinement $\widehat{\mu(A)}$ satisfies the twisted Maurer-Cartan equation
\(\label{twist}
\widehat{\mu(A)}=d\hat{A}- \overline{\hat{A}}\cdot {\hat{A}}\equiv B ~~{\rm mod} \ker(\hat{A})\;,
\)
where $B$ is some matrix in the ideal $M_{\rm form}$.
\end{theorem}
\theproof
Since $A$ is a formal connection, $\mu(A)$ satisfies
$$
\mu(A)=dA-\overline{A}\cdot A\equiv 0 ~~{\rm mod} \ker(A)\;.
$$
Hence, any refinement must satisfy 
$$
\widehat{\mu(A)}=D\hat{A}-\overline{\hat{A}}\cdot \hat{A}\equiv 0 ~~{\rm mod}\ \widehat{\ker(A)}\;,
$$
where $D=d+(-1)^*\delta$ is the {\v C}ech-Deligne differential on $M_{\rm diff}$. 
Now by sequence \eqref{ses}, we see that this is equivalent to existence of a matrix of forms $B$ satisfying \eqref{twist}.
\endofproof

In general, the Deligne-Beilinson cup product does not refine the de Rham wedge product 
for the whole de Rham complex, but does so only for the top and bottom degrees, as we have
seen in Prop. \ref{BD refines} and Prop. \ref{DB refines2}. However, for the triple product  the only cup products that arise
are between degree zero and degree one cocycles, so that nothing is missed in passing to $\cup_{DB}$. 
Consequently, for the case of the triple product, the matrix $B$ in the above example encodes the information needed to define the de Rham Massey product. 
More precisely, we have the following.

\begin{proposition} 
\label{Prop Msing}
Let $a_{i}\in H^{*}(X,\ZZ)$, $i=1,2,3$, and let $\iota(a)_i\in H_{\rm dR}^{*}(X)$ denote the inclusions into de Rham cohomology. Let
\begin{displaymath}
A= \left( \begin{array}{cccc}
0 & a_1 & \phi_{1,2} & \ast \\
0 & 0 & a_2 & \phi_{2,3} \\
0 & 0 & 0 & a_3 \\
0 & 0 & 0 & 0
\end{array} \right)\;.
\end{displaymath} 
in $M_{\rm sing}$ be a matrix of singular cochains defining a formal connection and let $\mu(A)$ be the corresponding solution to the corresponding Maurer-Cartan equation. Then for any differential refinement $\widehat{\mu(A)}$ of $\mu(A)$, the curvature $R(\widehat{\mu(A)})$ is a \emph{de Rham} Massey product in $\langle \iota(a)_1,\iota(a)_2,\iota(a)_3\rangle$. If, in addition, $\widehat{\mu(A)}$ is a solution to the \emph{differential} Maurer-Cartan equation, then $R(\widehat{\mu(A)})=0$ and $\mu(A)$ represents a torsion class.
\end{proposition} 
\theproof
Let $a_{i}$, $i=1,2,3$, be singular cochains of degree $n_{i,i+1}$. Suppose 
that the triple product is defined, and choose a defining system
\begin{displaymath}
A= \left( \begin{array}{cccc}
0 & a_1 & \phi_{1,2} & \ast \\
0 & 0 & a_2 & \phi_{2,3} \\
0 & 0 & 0 & a_3 \\
0 & 0 & 0 & 0
\end{array} \right)\;.
\end{displaymath}
Let 
\begin{displaymath}
\hat{A}= \left( \begin{array}{cccc}
0 & \hat{a}_1 & \hat{\phi}_{1,2} & \ast \\
0 & 0 & \hat{a}_2 & \hat{\phi}_{2,3} \\
0 & 0 & 0 & \hat{a}_3 \\
0 & 0 & 0 & 0
\end{array} \right)
\end{displaymath}
be a refinement. Then we know that the refinement $\widehat{\mu(A)}$ satisfies the equation
$D\hat{A}=\hat{A}\cdot \overline{\hat{A}}+B$
up to some element in $\ker(\hat{A})$. Explicitly, letting $B=(\eta_{ij})$, we have
\begin{displaymath}
\hat{A}= \left( \begin{array}{cccc}
0 & 0 & D\hat{\phi}_{1,2} & \ast \\
0 & 0 & 0 & D\hat{\phi}_{2,3} \\
0 & 0 & 0 & 0 \\
0 & 0 & 0 & 0
\end{array} \right)
= \left( \begin{array}{cccc}
0 & \eta_{12} & \hat{a}_1\cup \hat{a}_2+\hat{\eta}_{13} & \hat{a}_1\cup\hat{\phi}_{2,3}- \hat{\phi}_{1,2}\cup \hat{a}_3 \\
0 & 0 & \eta_{23} & \hat{a}_1\cup \hat{a}_2+\hat{\eta}_{24}  \\
0 & 0 & 0 & \eta_{34} \\
0 & 0 & 0 & 0
\end{array} \right)\;.
\end{displaymath}
The requirement that this equation hold up to an element in $\ker(\hat{A})$ forces the equations
\begin{eqnarray*} 
 \eta_{12}&=&0\;, \qquad D\hat{\phi}_{1,2}=\hat{a}_1\cup \hat{a}_2+\eta_{13}\;,
 \\
 \eta_{23}&=&0\;, \qquad D\hat{\phi}_{2,3}=\hat{a}_1\cup \hat{a}_2+\eta_{24}\;.
\end{eqnarray*}
At the level of connections, the data provided by the right two equations reduces to
\begin{eqnarray}\label{masseyeq}
d\phi_{1,2} &=&b_1\wedge a_2+\eta_{13}
\\
\label{masseyeq2}
d\phi_{2,3} &=& b_2\wedge a_3+\eta_{24}\;,
\end{eqnarray}
where $b_1$ and $b_2$ are forms representing the connections with curvatures $a_1$ and $a_2$.

Now forming $\widehat{\mu(A)}$ gives the matrix 
\begin{displaymath}
\widehat{\mu(A)}= \left( \begin{array}{cccc}
0 & 0 & 0 & \hat{a}_1\cup\hat{\phi}_{2,3}- \hat{\phi}_{1,2}\cup \hat{a}_3 \\
0 & 0 & 0 & 0  \\
0 & 0 & 0 & 0 \\
0 & 0 & 0 & 0
\end{array} \right)\;.
\end{displaymath}
Finally, applying the curvature map $R$ to the only nonzero term gives
$$
\hspace{-1mm}
\begin{aligned}
R\left(\hat{a}_1\cup\hat{\phi}_{2,3}-\hat{\phi}_{1,2}\cup \hat{a}_3\right) &= R(\hat{a}_1)\wedge R(\hat{\phi}_{2,3})-R(\hat{\phi_{1,2}})\wedge R(\hat{a}_3)
\\
&= R(\hat{a}_1)\wedge R(\hat{a}_2\cup \hat{a}_3)+R(\hat{a}_1)\wedge \eta_{24}-\left(\eta_{13}\wedge R(\hat{a}_3)+R(\hat{a}_1\cup \hat{a}_2)\wedge R(\hat{a}_3)\right)
\\
&= R(\hat{a}_1\cup (\hat{a}_2\cup \hat{a}_3))+R(\hat{a}_1)\wedge \eta_{24}-\eta_{13}\wedge R(\hat{a}_3)-R((\hat{a}_1\cup \hat{a}_2)\cup\hat{a}_3)
\\
&= R(\hat{a}_1)\wedge \eta_{24}-\eta_{13}\wedge R(\hat{a}_3)
\\
&= {a}_1\wedge \eta_{24}-\eta_{13}\wedge {a}_3\;.
\end{aligned}
$$
Notice that it follows from equations \eqref{masseyeq} and \eqref{masseyeq2} that the last line represents a de Rham Massey product (simply apply $d$ to both sides of those equations). This proves the first claim. 

For the second, observe that if $\hat{\mu}$ solves the Maurer-Cartan equation, then we can choose $B=(\eta_{ij})=0$,
 and the curvature calculated above must vanish.  
\endofproof

\section{Applications}
\label{Sec applications}


We will discuss our applications in this section, both from geometry and mathematical physics. 
We will show how Massey products arise in various settings, both classically and then in the 
newly constructed stacky form.

\subsection{Trivializations for (higher) structures}  
\label{Sec triv}

In this section we will consider Massey products arising from characteristic 
classes, hence associated with bundles or (higher) abelian gerbes. The refined 
Massey products will be associated with bundles or (higher) abelian gerbes 
together with connections on them. 
We consider examples involving the Deligne derivative $D$, which in the setting 
of the \v{C}ech-Deligne double 
complex, is given by  $D=d + (-1)^k \delta$.

\begin{example}
\label{Ex firstc}
Let $\pi:E\to M$ be a vector bundle equipped with connection $\nabla$. Let $\hat{c}_1(E,\nabla)$ be the {\v C}ech-Deligne cochain representing the differential refinement of the charateristic form corresponding to the connection (see \cite{Bun1}). Suppose that $\hat{c}_1(E,\nabla)$ is trivializable as a {\v C}ech-Deligne cochain and that moreover that there are cochains $\hat{a}$ and $\hat{b}$ such that $\hat{c}_{1}(E,\nabla)=\hat{a}\cup\hat{b}$. Since the class of $\hat{c}_1(E)$ vanishes in differential cohomology, there is a {\v C}ech-Deligne cocycle $\hat{A}$, with curvature $A$, such that 
\(\label{massey1}
D \hat{A}=\hat{c}_1(E,\nabla)=\hat{a}\cup\hat{b}.
\)

It was shown by Gomi \cite{Gom} that for a differential cohomology 
classes $\hat{a}$ of odd degree $n$, we have the formula
\(\label{steenrod trivialize}
[\hat{a}\cup\hat{a}]=ji(Sq^{n-1}(\bar{a}))\;.
\)
Here, $i$ is the map on cohomology induced via the representation as the square roots of unity $i:\ZZ/2\to U(1)$ (see \cite{Gom}\cite{Bun1}), and $j$ denotes the inclusion into differential cohomology via the 
map in diagram \eqref{diamond} which raises the degree by 1.  Let us assume that $a$ is divisible by 2 so that the 
mod 2-reduction is trivial and choose a trivializing {\v C}ech cochain $\phi$. Write $\varphi=ji(\phi)$. In this case, \eqref{steenrod trivialize} implies the equation
\(\label{massey2}
D(\varphi)=ji\delta(\phi)=\hat{a}\cup\hat{a}\;.
\)
Now the following matrix organizes the defining system given by 
equations \eqref{massey1} and \eqref{massey2}:
\begin{displaymath}
\left( \begin{array}{cccc}
0 & \hat{a} & \varphi &  \\
0 & 0 & \hat{a} & \hat{A} \\
0 & 0 & 0 & \hat{b} \\
0 & 0 & 0 & 0
\end{array} \right)\;.
\end{displaymath}
Then an element of the Massey product $\langle \hat{a},\hat{a},\hat{b}\rangle$ is given by the class of the
 {\v C}ech-Deligne cochain
$$
\hat{a}\cup\hat{A}-\varphi\cup \hat{b}\;,
$$
which is an element in $\widehat{H}^2(E;\ZZ)$. 
\end{example}


The previous example can be generalized to higher Chern classes. 
%

\begin{example}
Let $E\to M$ be a vector bundle with connection $\nabla$. Suppose that at the level of {\v C}ech-Deligne cochains, we have
\( 
\hat{c}_{2n-1}(E,\nabla)=\hat{a}_{2n-1}\cup\hat{b}_{2n-1}  
\qquad 
\text{and}
\qquad 
D\hat{A}_{(4n-3)}=\hat{c}_{2n-1}(E,\nabla)=\hat{a}_{2n-1}\cup\hat{b}_{2n-1}\;,
\label{eq cn}
\)
so that $\hat{c}_{2n-1}(E,\nabla)$ is trivializable as a bundle equipped with connection. We also assume that
$Sq^{n-1}(\bar{a}_{2n-1})=0$,
where $\bar{a}$ is the mod 2 reduction of $a$. Then, as in example \ref{Ex firstc} we have 
$\hat{a}_{2n-1}\cup\hat{a}_{2n-1}=D\varphi$,
for some cochain $\varphi$ \cite{Gom}. We have
\(\label{masseyn}
D(\varphi)=ji(Sq^{n-1}(\bar{a}_{2n-1}))=\hat{a}_{2n-1}\cup\hat{a}_{2n-1}.
\)
Now the following matrix organizes the defining system given by 
equations \eqref{eq cn} and \eqref{masseyn}:
\begin{displaymath}
\left( \begin{array}{cccc}
0 & \hat{a}_{2n-1} &\varphi &  \\
0 & 0 & \hat{a}_{2n-1} & \hat{A}_{(2n-1)} \\
0 & 0 & 0 & \hat{b}_{2n-1} \\
0 & 0 & 0 & 0
\end{array} \right)\;,
\end{displaymath}
and an element of the Massey product $\langle \hat{a}_{2n-1},\hat{a}_{2n-1},\hat{b}_{2n-1}\rangle$ is given by the class
$$\hat{a}_{2n-1}\cup\hat{A}_{(2n-1)}-\varphi\cup \hat{b}_{2n-1}\;.$$
\end{example}

We now consider the more interesting trivializations of String, Fivebrane \cite{SSS2}  
and Ninebrane structures \cite{9brane}. 
In fact, what we will consider are slightly weaker
versions, i.e. the vanishing of the $p_i$, $i=1,2,3$, where $p_i$ is the $i$th Pontrjagin class 
rather than the vanishing of the precise fractional classes. These differ from $p_i$-structures
by the fact that we still require the lower Pontrjagin classes to vanish
(see \cite{9brane} for more discussion). We will then in turn consider differential refinements 
of these structures, leading to Massey products representing \emph{geometric} String, Fivebrane  
and Ninebrane structures, respectively.

\begin{example}[Differential String structures and Chern-Simons theory]
\label{p1cs3}
On a smooth manifold $M$, viewed as a stack, consider a Spin bundle $E$ with 
connection $\nabla$  characterized by a morphism of stacks 
$\nabla: M \to \BB {\rm Spin}(n)_{\rm conn}$, 
to the moduli stack of bundles of rank $n$ Spin bundles with Spin connections. At the level of classifying 
spaces, the fractional Pontrjagin class appears as a map
$$
\tfrac{p_1}{2}:B{\rm Spin}(n)\longrightarrow B^3U(1)\simeq K(\ZZ,4)
$$
which obstructs String orientability. There is a unique differential refinement of the first Spin characteristic class 
$\tfrac{p_1}{2}$ denoted $\widehat{\tfrac{p_1}{2}}$ which gives a map at the level of moduli stacks
$$
\widehat{\tfrac{p_1}{2}}:\BB {\rm Spin}(n)_{\rm conn}\longrightarrow \BB^3U(1)_{\rm conn}\;.
$$
and captures the data of Chern-Simons theory (see \cite{FSSt} \cite{SSS3} \cite{FSS1} \cite{FSS2}
\cite{Bu5} \cite{Wal} \cite{Red} \cite{CJMSW}). 
Composing this map with with a map $\nabla:M\to \BB {\rm Spin}_{\rm conn}$ giving a Spin bundle, equipped with connection and resolving $M$ by its {\v C}ech nerve gives a {\v C}ech-Deligne cochain
$\widehat{\tfrac{p_1}{2}}(\nabla)$ on $M$. Suppose that the Spin bundle trivializes as a bundle with connection, i.e. 
that we have $\widehat{\tfrac{p_1}{2}}(\nabla)=0$ as a differential cohomology class. There are two interesting cases that can 
arise in practice and we will treat these cases separately. Suppose that $\widehat{\tfrac{p_1}{2}}(\nabla)$ decomposes as 
a square of a {\v C}ech-Deligne cochain. That is, we have 
\(
\widehat{\tfrac{p_1}{2}}(\nabla)=\hat{a} \cup_{\rm DB} \hat{a}\;.
\label{pDB}
\)
Diagrammatically, we have 
\(\label{trivialize p1}
\xymatrix{
&& 0 \ar@/^1pc/[drr]_<<<<{\ }="s" && 
\\
M \ar@/_1pc/[drr]_-{(\hat{a}, \hat{a})} \ar@/^1pc/[urr]
 \ar[rr]^-<<<<<<<<<<<<<<<<<<<{\nabla\ \ \ \ \ }="t" && \BB {\rm Spin}(n)_{\rm conn} \ar[rr]^-{\widehat{\tfrac{p_1}{2}}} &&
\BB^3U(1)_{\rm conn}\;,
\\
&& \BB U(1)_{\rm conn} \times \BB U(1)_{\rm conn} 
\ar@/_1pc/[urr]_-{\cup_{\rm DB}}
\ar@{=>} "s"; "t"
}
\)
where, by the trivialization condition \eqref{pDB}, the lower diagram commutes strictly,
and when we pass to connected components $\pi_0{\rm Map}(M, \BB^3U(1)_{\rm conn})$
the map ${\widehat{\tfrac{p_1}{2}}}$ is trivial, so that the upper part of the 
diagram commutes up to homotopy. A choice of homotopy is precisely a
trivializing {\v C}ech-Deligne 3-cochain $\hat{B}$. Given two such cochains $\hat{B}$ and $\hat{C}$, the difference is necessarily a cocycle since
$$D(\hat{B}-\hat{C})={\widehat{\tfrac{p_1}{2}}}-{\widehat{\tfrac{p_1}{2}}}=0\;.$$
Consider the defining system
$$
\left( \begin{array}{cccc}
0 & \hat{a} &\hat{B} &  \\
0 & 0 & \hat{a} & \hat{C} \\
0 & 0 & 0 & \hat{a} \\
0 & 0 & 0 & 0
\end{array} \right)\;.
$$
The corresponding Massey product then takes the form
\(
\langle \hat{a}, \hat{a}, \hat{a} \rangle = \hat{B} \cup_{\rm DB} \hat{a} 
- \hat{a} \cup_{\rm DB} \hat{C}=\hat{a}\cup_{\rm DB}(\hat{B}-\hat{C})\;.
\label{MaCS}
\)
Thus we can identify the Massey product as a flat bundle which is built entirely out of the trivializations of the ${\rm Spin}$ bundle with connection $\nabla$. Another interesting case happens when ${\widehat{\tfrac{p_1}{2}}}$ decomposes as $\hat{a}\cup_{DB}\hat{b}$. In this case, if the class of both $\hat{a}\cup_{DB}\hat{a}$ and ${\widehat{\tfrac{p_1}{2}}}$ vanish in differential cohomology, choosing local trivialization $\hat{B}$ and $\hat{C}$ of $\hat{a}\cup_{\rm DB}\hat{a}$ and ${\widehat{\tfrac{p_1}{2}}}$ (respectively) lead to the defining system
$$
\left( \begin{array}{cccc}
0 & \hat{a} &\hat{B} &  \\
0 & 0 & \hat{a} & \hat{C} \\
0 & 0 & 0 & \hat{b} \\
0 & 0 & 0 & 0
\end{array} \right)\;,
$$
and we get the Massey product
$$
\langle \hat{a}, \hat{a}, \hat{b} \rangle = \hat{B} \cup_{\rm DB} \hat{b} 
- \hat{a} \cup_{\rm DB} \hat{C}\;.
\label{MaCS}
$$
In this case the trivialization of the ${\rm Spin}$ bundle and the trivialization of the square $\hat{a}\cup_{\rm DB}\hat{a}$ combine to give a flat bundle representing the Massey product.
\end{example}

\begin{remark} {\bf (i)} Note that the above example can be extended to 
the case when the Spin bundle has a different rank than the dimension of 
the manifold. In particular, this holds for the stable case. 

\item {\bf (ii)} Note that \eqref{pDB} 
implies, in particular, that at the level of de Rham cohomology 
we locally have  
$dB_2=CS_3(\nabla)$, where
$B_2$ is the connection on the bundle ${\hat B}$. This then can be viewed
as a generalization of local trivialization of Chern-Simons theory. Hence the 
Massey product is s bundle on $E$ that is built out of the trivializations, 
including those of Chern-Simons. 
Furthermore, the structure of the Massey product \eqref{MaCS} indicates that,
even though we have a trivialization of Chern-Simon theory, we still have 
some secondary structure. 
\footnote{Note that Chern-Simons theory by itself can be viewed in a sense
as a secondary structure, so the above is a secondary structure (in one
sense) on some other secondary structure. We plan to make this precise elsewhere.} 

\item {\bf (iii)} Note that Example \ref{p1cs3} generalizes in a similar fashion to 
the cases of differential Fivebrane \cite{SSS3} and differential Ninebrane structures
\cite{9brane} with trivializing conditions on the 
characteristic classes given by $\hat{\tfrac{p_2}{6}}(\nabla)=D\hat{B}_6=D\hat{C}_6$ and 
$\hat{\tfrac{p_3}{240}}( \nabla)=D\hat{B}_{10}=D\hat{C}_10$,
respectively, with trivializing bundles $\hat{B}_i\neq \hat{C}_i$ of degree $i$. If $\hat{\tfrac{p_3}{240}}( \nabla)$ decomposes as the square $\hat{a}\cup_{\rm DB}\hat{a}$, the diagram \eqref{trivialize p1} 
will have the obvious modifications in degrees with the middle entry being replaced
by the appropriate structure, e.g. $\BB {\rm String}_{\rm conn}$ for the case of a
Fivebrane structure. 
The trivialization of these structures a priori give rise to Chern-Simons theories
in dimension 7 and 11, respectively, as highlighted in \cite{SSS2} \cite{9brane}. 
In the current setting, we will have trivializations of the Chern-Simons theories 
themselves at the level of complete data of bundles with connections, and governed 
by the corresponding Massey products, which would read the same as \eqref{MaCS} 
but with obvious changes in degrees. 
\end{remark}

\begin{remark} 
[Transfer of Massey products] {\bf (i)} 
A natural question is whether one can relate the stacky Massey triple product to the triple Deligne-Beilinson 
cup product. To that end, we recall the following from \cite{KS} (the argument there was for 
specific dimensions 
but it extends evidently to any dimension). 
Consider $Z^{n+1}$ as obtained
from gluing two cobordisms together, i.e. $Z^{n+1}$ is an
orientable compact manifold and $Y^{n}$ is a submanifold of
codimension 1 such that $Z^{n+1}-Y^n$ has two connected components, each
of which is a cobordism. Then from the Mayer-Vietoris sequence,
there is a connecting (or transfer) map 
\(
  T: H^k(Y^{n}) \longrightarrow H^{k+1}(Z^{n+1})\;.
\)  
Now let
$a,b,c\in H^*(Z^{n+1})$ with restrictions
$a^\prime, b^\prime, c^\prime \in H^*(Y^{n})$, 
and suppose
further that the cup products vanish $a^{\prime} \cup b^{\prime} =b^{\prime}\cup c^{\prime}=0\in H^*(Y^{n})$ so that the Massey product is defined. 
Then, by considering the Poincar\'e dual chains, one has that the transfer of the Massey product gives
the triple product \cite{KS}
 \(
 \label{emtt}
 {T\langle a^{\prime},b^{\prime},c^{\prime}\rangle= a\cup b\cup c \mod
\text{indeterminacy}}
\)
 where the Massey product is taken in $H^*(Y^{n})$,
the product in $H^*(Z^{n+1})$. The indeterminacy can be taken as $a\cup z+x\cup c$
where $z,x$ are cocycles in the opposite connected components of
$Z^{n+1}-Y^{n}$.
We propose generalizing this to our stacky setting of differential cohomology. 
We expect that the connecting homomorphism for differential cohomology takes the form
$$
  T: H^{k-1}(Y^{n},U(1)) \longrightarrow \hat{H}^{k+1}(Z^{n+1})\;,
$$
and sends the differential Massey product to the triple DB cup product (modulo indeterminacy).
%
\item {\bf (ii)} The Deligne-Beilinson triple cup product arises in the description of certain 
Chern-Simons type field theories in \cite{FSS1} \cite{FSS2}. The above then would be 
applied to these theories, giving that 
the Massey triple product of three differential cohomology elements on $Z^{n+1}$ 
transfers to a triple cup product Chern-Simons theory (in the sense of \cite{FSS1} \cite{FSS2})
on $Y^n$.
We leave the details of checking this for the future. 
\end{remark}

\subsection{Characteristic forms and anomaly cancellations} 
\label{Sec anomaly}

Presence of anomalies in a physical theory parametrizes to which extent 
certain entities are not (well) defined. Cancellation of these anomalies amounts to 
defining physical entities in the right mathematical setting. The process often requires 
an extension  of a topological or geometric setting to a more refined one. For example, 
to be able to talk about spinors, one has to set up the problem in the Spin bundle
as opposed to the tangent bundle. This requirement is obstructed by the second 
Stiefel-Whitney class, and the structure itself leads to interesting geometry and topology. 
One important instance of this is the Green-Schwarz anomaly cancellation condition 
required for consistency of string theory, which from the mathematical point of view
essentially requires working on manifolds with a (twisted) String structure. 
See \cite{Freed}  \cite{SSS2} \cite{Matt} for readable accounts aimed at mathematicians. 

\medskip
A generic situation is as follows.
Consider a bundle $P$ with curvature $F$ on a manifold $M$. Let $c_i(P)$ be a characteristic 
class of degree $i$ and let $c_i(F)$ be the corresponding characteristic form. Consider the 
conditions in cohomology $c_i(P) \cup c_j(P)=0$ and $c_j(P) \cup c_k(P)=0$. Then at the level 
of characteristic forms we have the trivializations via differential forms $\alpha$ and $\beta$ 
of the indicated degrees
\(
\label{trivchar}
c_i(F) \wedge c_j(F) =d\alpha_{(i+j-1)}\;,
\qquad \qquad
c_j(F) \wedge c_k(F) =d\beta_{(j+k-1)}\;.
\)
We build the composite differential form 
$$
\mu=c_i(F) \wedge \beta_{(j+k-1)} + (-1)^{i-1} \alpha_{(i+j-1)} \wedge c_k(F) ~ \in \Omega^{i+j+k-1}(M)\;,
$$
which is directly verified to be closed. This then allows us to form the Massey triple product 
of the corresponding cohomology classes
$$
\langle c_i(P), c_j(P), c_k(P) \rangle \in H^{i+j+k-1}(M; \Z)\;.
$$
Notice that we can consider conditions analogous to \eqref{trivchar} in differential cohomology
\(
\label{triv char}
\hat{c}_i(F) \cup \hat{c}_j(F) =D\hat{\alpha}_{(i+j-1)}\;,
\qquad \qquad
\hat{c}_j(F) \cup \hat{c}_k(F) =D\hat{\beta}_{(j+k-1)}\;,
\)
requiring not only that the characteristic forms vanish, but that the corresponding bundles 
trivialize as bundles with connection. In this case, we can form the bundle (differential cochain)
\(
\hat{\mu}=\hat{c}_i(F) \cup \hat{\beta}_{(j+k-1)} + (-1)^{i-1} \hat{\alpha}_{(i+j-1)} \cup \hat{c}_k(F) 
~ \in \map(M,\BB^{i+j+k-2}U(1)_{\rm conn})\;,
\label{eq mu hat}
\)
which is an element in
\(
\langle \hat{c}_i(P), \hat{c}_j(P), \hat{c}_k(P) \rangle \in \hat{H}^{i+j+k-1}(M; \Z)\;.
\)
We summarize the above.
\begin{proposition}
Given a system \eqref{triv char} of  trivializations of products of differential characteristic classes,  
we can build the stacky Massey product given by \eqref{eq mu hat}. 
\end{proposition}

We now provide an application of this direct but fairly general observation.
Consider a 10-dimensional manifold $X^{10}$ with metric $g$ on which there is a vector
bundle with connection $A$. One can consider the setting in families, i.e. take a bundle 
$E$ with fiber $X^{10}$ and base a parameter space and then integrate over the fiber
to get a class on the parameter space (see \cite{Freed} for beautiful 
constructions). We will not do all this but simply just set up integral expressions 
which will suffice for our purposes. 
The Green-Schwarz anomaly polynomials are given as
\bea
I_4 &=& p_1(g) - {\rm ch}_2(A)\;,
\nonumber\\
I_8 &=& -{\rm ch}_4(A) + \tfrac{1}{48} p_1(g) {\rm ch}_2(A) -
\tfrac{1}{64} p_1(g)^2 + \tfrac{1}{48}p_2(g)\;.
\eea
In \cite{SSS3} the first polynomial $I_4$ is interpreted as giving rise to a twisted String structure, and the 
indecomposable terms $p_2(g)$ and ${\rm ch}_2(A)$ in $I_8$ are 
interpreted as giving rise (essentially) to a Fivebrane structure and its twist, respectively. 
Their trivializations $H_3$ and $H_7$ provide trivializations of String and 
Fivebrane structures, respectively. 
A question remained on how to interpret the decomposable terms in $I_8$, namely 
 $\tfrac{1}{48} p_1(g) {\rm ch}_2(A)$ and $-\tfrac{1}{64} p_1(g)^2$. 
 We provide one interpretation of the corresponding trivializations, which fits well 
within our context. Consider the situation when $[p_1\cup {\rm ch}_2]=0=[p_1 \cup p_1]$,
i.e. 
\begin{eqnarray}
{\rm ch}_2(A) \wedge p_1(g) = d\alpha_7(A, g)\;,
\qquad \qquad 
p_1(g) \wedge p_1(g) = d\beta_7(g)\;,
\label{eq alpha}
\end{eqnarray}
and build the differential form 
\(\label{green}
\mu_{11}={\rm ch}_2(A) \wedge \beta_7(g) - \alpha_7(A, g) \wedge p_1(g)\;.
\)
This form is closed by virtue of \eqref{eq alpha}. Therefore, 
we can form the Massey triple product 
\(
\langle {\rm ch}_2, p_1, p_1 \rangle \in H^{11}(E; \Z)\;.
\)
As expected, the previous discussion refines to differential cohomology. Let $X^{10}$ be as before. 
Since we are fixing a Riemannian metric on $X^{10}$ and equipping the vector bundle with a connection 
$A$, it follows by uniqueness of characteristic forms (see \cite{Bun1} \cite{SS})
that we have unique differential refinements
\begin{eqnarray}
\hat{I}_4 &=& \hat{p}_1(g) - {\rm \hat{ch}}_2(A)\;,
\label{mixed first}
\\
\hat{I}_8 &=& -{\rm \hat{ch}}_4(A) + \tfrac{1}{48} \hat{p}_1(g) {\rm \hat{ch}}_2(A) -
\tfrac{1}{64} \hat{p}_1(g)^2 + \tfrac{1}{48}\hat{p}_2(g)\;.
\label{mixed I4I8}
\end{eqnarray}
We now consider the situation when these bundles trivialize as bundles with connections:
$[\hat{p}_1\cup {\rm \hat{ch}}_2]=0=[\hat{p}_1 \cup \hat{p}_1]$, so that expressions \eqref{eq alpha}
get replaced by
$$
{\rm \hat{ch}}_2(A) \wedge \hat{p}_1(g) = D\hat{\alpha}_7(A, g)\;,
\qquad \qquad
\hat{p}_1(g) \wedge \hat{p}_1(g) = D\hat{\beta}_7(g)\;.
$$
We then build the bundle
\(
\hat{\mu}_{11}={\rm \hat{ch}}_2(A) \cup \hat{\beta}_7(g) - \hat{\alpha}_7(A, g) 
\cup \hat{p}_1(g) ~\in \map(X^{10},\BB^{10}U(1)_{\rm conn})\;,
\)
which is a representative of the Massey triple product
\(
\langle {\rm \hat{ch}}_2, \hat{p}_1, \hat{p}_1 \rangle \in \hat{H}^{11}(E; \Z)\;.
\label{GS Massey}
\)
Therefore, we have the following
\begin{proposition}
The mixed terms in the Green-Schwarz anomaly polynomials \eqref{mixed first} \eqref{mixed I4I8} 
give rise to a stacky Massey product  given by the top class \eqref{GS Massey}. 
\end{proposition} 
It is interesting to note the form of the connection on the bundle $\hat{\mu}$. 
Using the formula for the DB cup product, we see that the connection is
\(\label{connect}
CS_3(A) \wedge CS_3(g)\wedge p_1(g) - \alpha_7(A, g) \cup p_1(g)\;,
\)
which we will make use of below (see Prop. \ref{Prop connCS}).

\paragraph{Fiber integration of Massey products and anomaly line bundles}
In \cite{FSS1} \cite{FSS2}, a fiber integration map was defined by taking the usual fiber integration in cohomology, 
lifting to differential cohomology and then lifting to the internal hom in sheaves of positively graded chain 
complexes to produce a map 
$$
\int_{\Sigma^{k}}-:[N(C(\{U_i\}),\ZZ_{\mathcal{D}}^{\infty}[n]]\longrightarrow \ZZ^{\infty}_{\mathcal{D}}[n-k]\;.
$$
Here $\Sigma^{k}$ is a paracompact manifold of dimension $k$ and $C(\{U_i\})$ is the {\v C}ech nerve corresponding to a good open cover of $\Sigma^k$. The lifts are provided by the construction of Gomi and Terashima in \cite{GT}. Post-composing with the quasi-isomorphism provided by the exponential and applying the Dold-Kan functor gives a morphism of stacks in the form of holonomy 
\(
{\rm hol}_{\Sigma^{k}}:=\exp\left(2\pi i\int_{\Sigma_{k}}-\right):[\Sigma^{k},\BB^{n}U(1)_{\rm conn}]
\longrightarrow
 \BB^{n-k}U(1)_{\rm conn}\;.
\)  
Again in \cite{FSS1} \cite{FSS2}, it was observed that the abelian Chern-Simons action functional can be described by 
post-composing the cup product morphism with this holonomy map. In particular, for a manifold  $\Sigma^{4k+3}$, this composite induces an intersection pairing on differential cohomology
\(
(\hat{x},\hat{y})
\longrightarrow
 \exp\left(2\pi i \int_{\Sigma^{4k+3}}\hat{x}\cup\hat{y}\right)\;.
\)
For $k=0$ and $\hat{y}=\hat{x}$, this pairing gives the usual Chern-Simons action. We now would like
 to describe how to lift this morphism to the Massey product (when defined). 
 In fact, when the differential Massey product is defined, 
 we have a map 
\(
\langle \hat{x},\hat{y},\hat{z}\rangle_U~:~
\Sigma^k \times U \longrightarrow  \BB^{n_1+n_2+n_3+2}U(1)_{\rm conn}\;,
\label{Diff MP}
\)
which is natural in any test space $U$.
Hence, we can apply the fiber integration map. Since Massey products necessarily define \emph{flat} bundles, we see immediately that we have the following.
\begin{proposition} 
The integration over the fiber of the differential Massey product \eqref{Diff MP} can be identifies with a map
$$
e^{\left(2\pi i \int_{\Sigma^k}\langle\hat{x},\hat{y},\hat{z}\rangle\right)}{}_U~:~
U \longrightarrow
\BB^{n_1+n_2+n_3+2-k}U(1)_{\rm conn}\;,
$$
which is natural in $U$. Moreover, this map defines a flat bundle on $U$, and the map factors through the inclusion 
$$j:\flat \BB^{n_1+n_2+n_3+2-k}U(1)\into \BB^{n_1+n_2+n_3+2-k}U(1)_{\rm conn}\;.$$ 
\end{proposition}

\medskip
\begin{remark} {\bf (i)} The above construction can be generalized to higher Massey products, 
as we can fiber integrate any differential cohomology class of any degree, 
including those that are Massey products. 

\item {\bf (ii)} The notation $e^{\left(2\pi i \int_{\Sigma^k}\langle\hat{x},\hat{y},\hat{z}\rangle\right)}$ is slightly abusive, 
since this map may not be well-defined on the entire Massey product (due to indeterminacy). What we really mean here is an \emph{element} of the Massey product. 
\end{remark}

In particular, when $\hat{x}$, $\hat{y}$ and $\hat{z}$ come as characteristic forms, they are given by \emph{morphisms} of stacks;
e.g.
$$
\hat{x}:[\Sigma^k,\BB G_{\rm conn}]\longrightarrow
 [\Sigma^k, \BB^{n_1}U(1)_{\rm conn}]\;,
$$
which gives a natural assignment of differential cohomology classes as we vary the $G$-principal 
bundle with connection on $\Sigma^k$ \cite{FSS1}. In this case, after choosing trivialization of 
$\hat{x} \cup_{\rm DB} \hat{y}$ 
and $\hat{y} \cup_{\rm DB} \hat{z}$, 
fiber integration gives
the morphism
 of stacks 
 \(
e^{\left(2\pi i \int_{\Sigma^k}: \langle\hat{x},\hat{y},\hat{z}\rangle\right)}:~[\Sigma^k,\BB G_{\rm conn}]
\longrightarrow
 \flat \BB^{n_1+n_2+n_3+2-k}U(1)\;.
\label{mor BG}
\)
One interesting instance of this morphism comes from the previous example of Green-Schwarz anomaly polynomials. 
That is, we are interested in the triple product $\langle {\rm \hat{ch}}_2, \hat{p}_1, \hat{p}_1 \rangle$. 
In this case, we get a morphism
\(
e^{\left(2\pi i \int_{X^{10}}\langle{\rm \hat{ch}}_2,\hat{p}_1,\hat{p}_1\rangle\right)}:~
[X^{10},\BB G_{\rm conn}]
\longrightarrow
\flat \BB^{10-10}U(1)=U(1)^\delta\;,
\label{Massey U(1)}
\)
from the moduli stack of bundles on $X^{10}$ equipped with connection to smooth $U(1)$-valued functions. 
It is useful to unwind this map at the level of connections. Indeed, noting \eqref{connect}, 
we have at that level:

\begin{proposition} 
\label{Prop connCS}
The connection on the bundle prescribed by
 \eqref{mor BG} is given by the form
$$
\int_{X^{10}}CS_3(A) \wedge CS_3(g)\wedge p_1(g) - \alpha_7(A, g) \cup p_1(g)\;.
$$
\end{proposition} 

\medskip
\begin{remark}
 {\bf (i)} The exponential of the functional on the right, being built out of
Chern-Simons forms, 
 is indeed in $U(1)$.  
\item {\bf (ii)} As the structure of the functional in the proposition involves a product of two 
Chern-Simons forms, this suggests a formulation where $X^{10}$ is viewed
as a manifold of corners of codimension two, in the sense of the setting 
 in \cite{corners1} \cite{corners2}. We hope to take up this point of view elsewhere. 
\end{remark}

\subsection{Twisted cohomology and twisted Bianchi identities} 
\label{Sec Bianchi}

We consider the Ramond-Ramond (RR) fields in type IIA string theory on a ten-dimensional manifold  $X^{10}$ with a B-field,
whose curvature is a closed three-form $H_3$.  The RR fields of various degrees can be combined into the expression ${\cal F}=\sum_{i=0}^{5} F_{2i}$,
and satisfy the twisted Bianchi $d{\cal F}_n + H_3 \wedge {\cal F}_{n-2}=0$. In components,
\begin{eqnarray}
H_3 \wedge {F}_0&=&- d{F}_2\;,
\qquad \quad
H_3 \wedge {F}_2=- d{F}_4\;,
\qquad \quad
H_3 \wedge {F}_4=- d{F}_6\;,
\nonumber \\
H_3 \wedge {F}_6&=&- d{F}_8\;,
\qquad \quad
H_3 \wedge {F}_8=- d{F}_{10}\;,
\qquad  \quad
d{F}_0=0=d{F}_{10}\;.
\label{eq RR}
\end{eqnarray}

\begin{remark}
From these we will build expressions of degree ten.  
\item {\bf (i)} Considering the first and fifth expressions in \eqref{eq RR}, we can set up
the top differential form 
$$
\mu = F_0 \wedge F_{10} + F_2 \wedge F_8\;.
$$
This is closed by dimension reasons, so that we can form the triple Massey product 
$$
\langle  F_0, H_3, F_8 \rangle \in H^{10}(X^{10}; \Z)\;. 
$$
\item {\bf (ii)} Considering the second and fourth expressions in \eqref{eq RR}, we build the top form
$$
\mu'= F_2 \wedge F_8 + F_4 \wedge F_6\;.
$$
This is closed again by dimension reasons, and we can build the triple Massey product 
$$
\langle  F_2, H_3, F_6 \rangle \in H^{10}(X^{10}; \Z)\;. 
$$
\end{remark}

\medskip
We now would like to refine the previous discussion to differential cohomology. Notice that since $dF_{2i}\neq 0$, we cannot simply put hats everywhere and expect the equations to hold at the level of ordinary differential cohomology. Consequently, there are two directions we can go. First, we could try to form Massey products in \emph{twisted} differential cohomology, which is outside the scope of the present paper. Second, we can view the $F_{2i}$'s as 
\emph{improved} gauge invariant field strengths corresponding to potentials $C_{2i-1}$ with curvatures $G_{2i}$, which are not gauge invariant. We will expand on this latter point of view. To this end, we require that the potentials $C_{2i-1}$ satisfy
\(
dC_{n}+H_{3}\wedge C_{n-2}=0
\)
Notice that this equation implies that the improved field strengths $F_{2i}$ vanish, by definition. Combining the potentials into the single potential ${\cal C}=\sum_{i=0}^{3} C_{2i-1}$ we have, by assumption, the equations
\(
H_3 \wedge {C}_1=- d{C}_3\;, 
\qquad 
H_3 \wedge {C}_3=- d{C}_5\;,
\qquad
H_3 \wedge {C}_5=- d{C}_7\;.
\label{RR2conn}
\)
These equations can be viewed as conditions on the connections for differential refinements of the field strengths $G_{2i}$. Indeed, the full differentially refined equations read
\(
\hat{H}_3 \cup \hat{G}_2=- D\hat{G}_4\;,
\qquad
\hat{H}_3 \cup \hat{G}_4=- D\hat{G}_6\;,
\qquad
\hat{H}_3 \cup \hat{G}_6=- D\hat{G}_8\;.
\label{RR2}
\)
\begin{enumerate}[label= {\bf (\roman*)}]
\item 
Considering the first and third expressions in \eqref{RR2}, we can form
the bundle
$$
\hat{\mu} = \hat{G}_2 \cup \hat{G}_{8} + \hat{G}_4 \cup \hat{G}_6\;,
$$
with higher connection
\(
C_1\wedge G_8+C_3\wedge G_6=C_1\wedge H_3\wedge C_5+C_3\wedge H_3\wedge C_3\;.
\label{RR con 1}
\)
This bundle is an element in  the stacky Massey triple product 
\(
\langle  \hat{G}_2, \hat{H}_3, \hat{G}_6 \rangle \in \hat{H}^{10}(X^{10}; \Z)\;. 
\label{RR Mas 1}
\)
\item  Considering instead the first and second expressions in \eqref{RR2}, we form the higher bundle 
$$
\hat{\mu}'= \hat{G}_4 \cup \hat{G}_6 + \hat{G}_4\cup \hat{G}_4\;.
$$
with higher connection
\(
C_3\wedge G_6+C_5\wedge G_4=C_3\wedge H_3\wedge C_3+C_5\wedge H_3\wedge C_1\;.
\label{RR con 2}
\)
This is an element in the stacky triple Massey product 
\(
\langle  \hat{G}_4, \hat{H}_3, \hat{G}_4 \rangle \in \hat{H}^{10}(X^{10}; \Z)\;. 
\label{RR Mas 2}
\)
\end{enumerate}

\begin{proposition}
The system of twisted Bianchi identities for the differential RR fields leads to 
two higher bundles with connections \eqref{RR con 1} and \eqref{RR con 2} which are elements 
in the stacky Massey products in top degree \eqref{RR Mas 1} and \eqref{RR Mas 2}, respectively. 
\end{proposition}

It would be interesting to investigate the implications of these expressions to string theory. 
For now we just observe that, essentially and up to signs, $\mu_1$ and $\mu_2$ are part of the couplings that arise 
in calculating the topological partition function of the RR fields (in the case when $H_3=0$) \cite{DMW} \cite{BM}. 
While we do not pursue this here, we expect  $\hat{\mu}_1$ and $\hat{\mu_2}$ 
to be relevant for the calculation
of the partition function in the twisted differential case, $\hat{H}_3 \neq 0$, as well,
extending the twisted topological case in \cite{MoS} \cite{MaS}.

\subsection{Quadruple Massey products}
\label{Sec quad} 

We now consider a setting inspired by type IIB string theory. The main feature of this theory that concerns 
us here is that it has fields of odd degree, where the degree  three play a somewhat special role. 
Consider four fields as cohomology classes $h_3^{(i)}\in H^3(X; \Z)$, $i=1, \cdots, 4$, on a  ten-dimensional 
manifold $X^{10}$,
  and consider analogues of three composite (Ramond-Ramond) fields $F_5^{(j)}$, $j=1, 2,3$,
 such that 
 $$
 h_3^{(1)} \wedge h_3^{(2)}=-dF_5^{(3)}\;,
 \qquad \quad
  h_3^{(2)} \wedge h_3^{(3)}=-dF_5^{(1)}\;,
 \qquad \quad
  h_3^{(3)} \wedge h_3^{(4)}=-dF_5^{(2)}\;.
 $$
Then there are further composite (again analogues of Ramond-Ramond)  
fields $F_7^{(i)}$, $i=1, \cdots, 4$, such that
 \bea
 F_5^{(3)} \wedge h_3^{(3)}&=&-dF_7^{(3)}\;,
 \qquad \qquad 
  h_3^{(1)} \wedge F_5^{(1)}=-dF_7^{(1)}\;,
 \nonumber\\
  F_5^{(1)} \wedge h_3^{(4)}&=&-dF_7^{(4)}\;,
  \qquad \qquad
  h_3^{(2)} \wedge F_5^{(2)}=-dF_7^{(2)}\;.
 \eea
Then we will end up (see below) having the Massey quadruple product as the integer 
$$
\langle h_3^{(1)}, h_3^{(2)}, h_3^{(3)}, h_3^{(4)} \rangle:=
-F_7^{(3)} \wedge h_3^{(4)} - F_7^{(1)} \wedge h_3^{(4)} 
-F_5^{(3)} \wedge F_5^{(2)}
+ h_3^{(1)} \wedge F_7^{(4)} + h_3^{(1)} \wedge F_7^{(2)} \in H^{10}(X^{10}; \Z)\cong \Z\;.
$$

\vspace{3mm} 
We now elaborate on the above. 
We first start with the triple Massey product in the
current IIB string theory inspired context. Let $[h_3^{(i)}]\in H^3(X^{10})$ ($i=1,2,3$) be
non-zero cohomology classes such that $[h_3^{(1)}] \cup
[h_3^{(2)}]=0$ and $[h_3^{(2)}] \cup [h_3^{(3)}]=0$. For the cocycle
representatives $h_3^{(i)}$, write
\(
h_3^{(1)} \cup h_3^{(2)}=dF_5^{(1)}
\qquad \text{and} \qquad
h_3^{(2)} \cup h_3^{(3)}=dF_5^{(2)}\;. 
\label{1}
\)
Notice that from these two equations one gets immediately that
$d(F_5^{(1)}\cup h_3^{(3)} + h_3^{(1)} \cup F_5^{(2)})=0$ by a
straightforward application of the Leibnitz rule. This can then be
used to define the triple Massey product as the subset of $H^8(X^{10})$
given by
 \( \big\langle  [h_3^{(1)}], [h_3^{(2)}], [h_3^{(3)}]
\big\rangle =\left\{ \left[ F_5^{(1)}\cup h_3^{(3)}+h_3^{(1)}\cup
F_5^{(2)} \right] \right\}, \) where $h_3$ and $F_5$ run over all
possible choices above. The indeterminacy in the choice of the
representative $w=h_3^{(1)}\cup F_5^{(2)}+F_5^{(1)}\cup h_3^{(3)}$
for the triple product lies in the ideal
$\left([h_3^{(1)}],[h_3^{(2)}]\right)$.

\vspace{3mm}
In order to connect with the Massey 4-fold product, it
is good to rewrite the triple product in matrix form. The classes
$[h_3^{(i)}]$ and the elements $F_5^{(i)}$ can be encoded in a
matrix form

\( \left(\begin{array}{cccc}
0&a_{11}&a_{22}&\\
&0&a_{22}&a_{23}\\
&&0 &a_{33}\\
&&& 0
\end{array}\right)
=\left(
\begin{array}{cccc}
0&h_3^{(1)}&F_5^{(1)}& \\
 &0&h_3^{(2)}&F_5^{(2)}\\
 &&0 &h_3^{(3)}\\
&&& 0
\end{array}\right).
\)
The defining properties of the Massey triple product can be
expressed in a matrix multiplication as
\begin{eqnarray}
d\left(
\begin{array}{cccc}
0&h_3^{(1)}&F_5^{(1)}& \\
 &0&h_3^{(2)}&F_5^{(2)}\\
 &&0 &h_3^{(3)}\\
&&&0
\end{array}\right)&=&
\left(
\begin{array}{cccc}
0&0&h_3^{(1)}h_3^{(2)}& \\
 &0& 0& h_3^{(2)}h_3^{(3)}\\
 &&0 &0\\
&&& 0
\end{array}\right)\;,
\nonumber\\
\left(
\begin{array}{cccc}
0&h_3^{(1)}&F_5^{(1)}& \\
 &0&h_3^{(2)}&F_5^{(2)}\\
 &&0 &h_3^{(3)}\\
&&&0
\end{array}\right)^2
&=&
\left(
\begin{array}{cccc}
0&0&h_3^{(1)}h_3^{(2)}&h_3^{(1)}F_5^{(2)}+F_5^{(1)}h_3^{(3)} \\
 &0&0&h_3^{(2)}h_3^{(3)}\\
 & &0& 0
\end{array}\right)\;.
\end{eqnarray}

Now in order to go one step further to the quadruple
(or 4-fold) product, we need to satisfy certain conditions on the
triple product, in analogy to saying that higher obstructions arise
only once the lower ones vanish. So in our case, we first need to
assume that we can complete our set by adding two more elements, 
a fourth $h_3^{(4)}$ and a third $F_5^{(3)}$, such that 
\(
dF_5^{(3)}=h_3^{(3)}\cup h_3^{(4)}\;.
\label{3}
\) 
Besides the above representative $w$, we then have a second
representative for the triple product and is given by
$z=h_3^{(2)}\cup F_5^{(3)}+F_5^{(2)}\cup h_3^{(4)}$, namely
representing $\big\langle
[h_3^{(2)}],[h_3^{(3)}],[h_3^{(4)}]\big\rangle$. The condition to be
able to define the quadruple product is that {\it both} triple
products vanish {\it simultaneously}, i.e. that both cohomology
representatives $w$ and $z$ can be chosen as coboundaries, which we
write as $w=dF_7^{(1)}$ and $z=dF_7^{(2)}$.

\vspace{3mm}
We are now ready to define the 4-fold or quadruple
Massey product. In analogy to the triple product, we start with the
equations \eqref{1} and \eqref{3}, and then write the two
cocycles of degree eight
$$
dF_7^{(1)}=h_3^{(1)}\cup F_5^{(2)} + F_5^{(1)}\cup h_3^{(3)}
\qquad \text{and} \qquad 
dF_7^{(1)}=h_3^{(2)}\cup F_5^{(3)} + F_5^{(2)}\cup h_3^{(4)}\;,
$$
from which we get a cocycle \( x=h_3^{(1)}\cup
F_7^{(2)}+F_5^{(1)}\cup F_5^{(2)}+ F_7^{(1)}\cup h_3^{(4)} \) of
degree ten. 

\begin{remark}
{\bf (i)} Again, we define the quadruple Massey product $\big\langle
[h_3^{(1)}],[h_3^{(2)}],[h_3^{(3)}],[h_3^{(4)}]\big\rangle$ as a
collection of all cohomology classes $[x]\in H^{10}(X^{10})$ that we can
obtain by the above procedure.


\item {\bf (ii)} 
The indeterminacy is best presented as the {\it matrix} triple
product of certain elements,
%
namely of $(h_3^{(1)}, H^5(X^{10}))$,
$\begin{pmatrix}
h_3^{(2)} &
H^5(X^{10}) \\
0 & h_3^{(3)}
\end{pmatrix}$, 
and 
$\begin{pmatrix}
H^5(X^{10})\\
h_3^{(4)}
\end{pmatrix}$.
\end{remark}


\medskip
We now generalize this construction to differential cohomology. To produce the desired products, we again view the $F^{(i)}_5$ 
and $F^{(i)}_7$ as improved, gauge invariant field strengths and denote the corresponding potentials as $C^{(i)}_4$ and $C^{(i)}_6$, 
with curvatures $G^{(i)}_5$ and $G^{(i)}_7$. We now lift everything to the level of differential cohomology, which yields the equations
 \(
 \hat{h}_3^{(1)} \cup \hat{h}_3^{(2)}=-D\hat{G}_5^{(3)}\;,
 \qquad 
  \hat{h}_3^{(2)} \cup \hat{h}_3^{(3)}=-D\hat{G}_5^{(1)}\;,
 \qquad
  \hat{h}_3^{(3)} \cup \hat{h}_3^{(4)}=-D\hat{G}_5^{(2)}\;,
\label{first set}
\)
and 
 \begin{eqnarray}
 \hat{G}_5^{(3)} \cup \hat{h}_3^{(3)}&=&-D\hat{G}_7^{(3)}\;,
 \qquad \qquad
  \hat{h}_3^{(1)} \cup \hat{G}_5^{(1)}=-D\hat{G}_7^{(1)}\;,
 \nonumber\\
  \hat{G}_5^{(1)} \cup \hat{h}_3^{(4)}&=&-D\hat{G}_7^{(4)}\;,
  \qquad \qquad
  \hat{h}_3^{(2)} \cup \hat{G}_5^{(2)}=-D\hat{G}_7^{(2)}\;.
 \label{second set}
 \end{eqnarray} 
 The connection on the higher bundle is calculated as follows. Set
 $$
 \mathcal{A}:= -C_6^{(3)} \wedge h_3^{(4)} -  C_6^{(1)} \wedge h_3^{(4)} 
-C_4^{(3)} \wedge G_5^{(2)}
+ b_2^{(1)} \cup G_7^{(4)} + b_2^{(1)} \cup G_7^{(2)}
=
-C_6^{(3)} \wedge h_3^{(4)} - C_6^{(1)} \wedge h_3^{(4)}\;.
$$
Then, by writing the higher components in the last three terms 
via lower components, we get 
\(
\mathcal{A}=-C_6^{(3)} \wedge h_3^{(4)} - C_6^{(1)} \wedge h_3^{(4)} 
-C_4^{(3)} \wedge b^{(3)}_2\wedge h^{(4)}_3
+ b_2^{(1)} \wedge C^{(1)}_4\wedge h^{(4)}_3 +
 b_2^{(1)} \wedge b^{(2)}_2 \wedge b^{(3)}_2\wedge h^{(4)}_3\;.
\label{Higher conn}
\) 
Here $b_2^{(i)}$ denotes a local potentials for the forms $h^{(i)}_3$. 
 Therefore, we have the following description as phase or holonomy. 
 
 \begin{proposition}  
 The system \eqref{first set} \eqref{second set} 
 leads to the stacky Massey quadruple product 
$$
\langle \hat{h}_3^{(1)}, \hat{h}_3^{(2)}, \hat{h}_3^{(3)}, \hat{h}_3^{(4)} \rangle:=
-\hat{G}_7^{(3)} \cup \hat{h}_3^{(4)} - \hat{G}_7^{(1)} \cup \hat{h}_3^{(4)} 
-\hat{G}_5^{(3)} \cup \hat{G}_5^{(2)}
+ \hat{h}_3^{(1)} \cup \hat{G}_7^{(4)} + \hat{h}_3^{(1)} \cup \hat{G}_7^{(2)} 
\in \hat{H}^{10}(X^{10}; \Z)\;,
$$
viewed as a higher bundle whose connection ${\cal A}$ is given by \eqref{Higher conn}.
\end{proposition}

The discussion using matric Massey products carries over to differential cohomology in a similar fashion.
We also leave the discussion on the physical impact of the above constructions 
to a separate treatment. 

\vspace{6mm}
\noindent {\bf \large Acknowledgement}

\medskip
\noindent The authors would like to thank Domenico Fiorenza and Urs Schreiber for very useful 
discussions and comments, Chris Kapulkin for a useful comment
on the first version of the manuscript, and the referee for a careful reading of the manuscript and for many useful suggestions.


\end{document}